\documentclass[conf]{new-aiaa}
\usepackage[utf8]{inputenc}
\usepackage{graphicx}
\usepackage{amsmath, amsfonts}
\usepackage[nameinlink]{cleveref}
\usepackage{longtable, tabularx}
\usepackage{algorithm}
\usepackage[noend]{algpseudocode}
\usepackage{booktabs}
\usepackage{subcaption}
\usepackage[dvipsnames]{xcolor}
\usepackage{pgfplots}
\usepackage{tabto}
\usepackage{soul}
\usepackage{tikz}

\pgfplotsset{compat=1.18}
\usetikzlibrary{arrows.meta}

\definecolor{BYUblue}{RGB}{0, 61, 165}
\hypersetup{
    colorlinks=true,
    linkcolor=blue,
    citecolor=blue,
    urlcolor=blue,
    pdftitle={Fast Quadratic Manifold Learning For Nonlinear Dimensionality Reduction in Large-scale Systems using Riemannian Optimization},
    pdfauthor={Gavin Paxton, Seunghee Cheon, Rudy Geelen, Shane~A.~McQuarrie}
}

\usepackage{silence}
\definecolor{mplblue}{rgb}{0.1216, 0.4667, 0.7059}
\definecolor{mplorange}{rgb}{1.0, 0.498, 0.055}
\definecolor{mplgreen}{rgb}{0.173, 0.627, 0.173}

\algrenewcommand\algorithmicrequire{\textbf{Input:}}
\algrenewcommand\algorithmicensure{\textbf{Output:}}
\algnewcommand{\LineComment}[1]{\vspace{.5em}\noindent\hspace{-1.1em}\textcolor{gray}{\texttt{\# #1}}}
\algrenewcomment[1]{\tabto{1.5in}\textcolor{gray}{\texttt{\# #1}}}

\crefname{equation}{}{}
\Crefname{equation}{}{}
\crefname{figure}{Fig.}{}
\Crefname{figure}{Fig.}{}

\newcommand{\Vr}{\mathbf{V}_{\!r}}
\newcommand{\Vq}{\mathbf{V}_{\!q}}

\title{Fast Quadratic Manifold Learning For Nonlinear \\ Dimensionality Reduction in Large-scale Systems \\ using Riemannian Optimization}

\author{Gavin Paxton\footnote{Undergraduate researcher}, Seunghee Cheon\footnote{Graduate Student} and Rudy Geelen\footnote{AIAA Member}}
\affil{Department of Aerospace Engineering, Texas A\&M University, College Station, TX, 77843}

\author{Shane~A.~McQuarrie}
\affil{Department of Mathematics, Brigham Young University, Provo, UT, 84602}

\begin{document}

\maketitle

\begin{abstract}
The effectiveness of dimensionality reduction with quadratic manifolds hinges on the choice of a reduced basis and the associated quadratic correction terms. Existing approaches typically rely on subspaces spanned by the leading principal components of the training data. Although optimal for linear approximation, such bases are inherently suboptimal for quadratic manifold learning. Greedy basis-selection methods can significantly improve the representational capacity of quadratic manifolds by searching over a larger pool of candidate principal components, but the combinatorial cost limits the basis sizes that can be used in practice.
This work proposes \texttt{FastQM}, an approach that treats the identification of an optimal quadratic approximation as a continuous optimization problem on the Stiefel manifold. By rotating the reduced basis within a candidate span of singular vectors, \texttt{FastQM} learns an ideal coordinate alignment tailored to quadratic manifold approximation. A feature-space formulation ensures that the optimization cost scales independently of the full state-space dimension. The efficacy of the proposed method is demonstrated on a turbulent airfoil-wake large-eddy simulation.
\end{abstract}

\noindent

\section{Nomenclature}

{\renewcommand\arraystretch{1.0}
\noindent\begin{longtable*}{@{}l @{\quad=\quad} l@{}}
$t$  & Time variable \\
$N$ & Dimension of the high-dimensional state vector \\
$K$ & Number of training snapshots \\
$r$ & Dimension of the reduced-order linear subspace \\
$q$ & Dimension of the nonlinear correction subspace \\
$m$ & Number of candidate/feature modes \\
$\mathbf{s}$ & State vector, $\mathbf{s}(t) \in \mathbb{R}^N$\\
$\hat{\mathbf{s}}$ & Reduced-order state vector, $\hat{\mathbf{s}}(t) \in \mathbb{R}^r$ \\
$\mathbf{S}$ & Snapshot matrix, $\mathbf{S} \in \mathbb{R}^{N \times K}$ \\
$\tilde{\mathbf{S}}$ & Representation of the snapshot data in the candidate basis, $\tilde{\mathbf{S}} \in \mathbb{R}^{m \times K}$ \\
$\boldsymbol{\Xi}$ & Weight matrix for nonlinear correction terms, $\boldsymbol{\Xi} \in \mathbb{R}^{q \times r(r+1)/2}$ \\
$\gamma$ & Regularization parameter , $\gamma \geq 0 \in \mathbb{R} $ \\
$\sigma_j$ & $j$th singular value, $\sigma_i \geq 0 \in \mathbb{R}$ \\
$\Vr$ & Basis matrix for the linear term, $\Vr \in \mathbb{R}^{N \times r}$ \\
$\Vq$ & Basis matrix for quadratic correction term, $\Vq \in \mathbb{R}^{N \times q}$ \\
$\tilde{\mathbf{V}}$ & Basis matrix for candidate modes, $\tilde{\mathbf{V}} \in \mathbb{R}^{N \times m}$ \\
$\mathbf{Q}_r, \mathbf{Q}_q$ & Orthogonal matrices, $[~\mathbf{Q}_r~~\mathbf{Q}_q~] \in \text{St}(m, r+q)$ \\
$\mathbf{I}$ & Identity matrix \\
$\mathbf{0}$ & Zero matrix \\
$\operatorname{Tr}$ & Matrix trace operator \\ 
$\otimes$ & Kronecker product \\
$\|\cdot\|_F$ & Frobenius norm \\
$\operatorname{St}(m, r)$ & Stiefel manifold of $m \times r$ matrices with orthonormal columns \\
\end{longtable*}}

\section{Introduction}

\noindent

Large-scale physical systems are simulated by solving high-dimensional numerical models obtained by discretizing the governing equations, for example systems of partial differential equations (PDEs). Because such numerical models are often computationally expensive to solve, they are not typically suitable for real-time or multi-query applications such as design, optimization, or uncertainty quantification~\cite{benner2015promsurvey, peherstorfer2018survey}.
Model reduction addresses this challenge by constructing computationally efficient reduced-order models (ROMs) that approximate the dominant system behavior.
Because the fidelity of a ROM is fundamentally constrained by its underlying state-space representation, dimensionality reduction is a core pillar of ROM design.
Classical model reduction methods leverage linear subspaces to capture the essential dynamics; however, linear approximations often struggle to account for the intrinsic curvature of highly nonlinear solution manifolds. Consequently, building accurate ROMs may require large subspace dimensions, diminishing computational gains.
In this work, we introduce \texttt{FastQM}, a framework for learning \emph{quadratic} manifolds that minimizes reconstruction error by identifying a linear subspace and a matching quadratic correction term through numerical optimization.
Importantly, the dimensionality of the search space does not scale with the size of the original state space, ensuring that the computational cost remains manageable even for large-scale problems. Our numerical results demonstrate that \texttt{FastQM} outperforms proper orthogonal decomposition (POD) and POD-based quadratic manifolds in reconstruction accuracy, establishing it as a viable choice for large-scale model reduction tasks.

While POD~\cite{lumley1967structure, sirovich1987turbulence, berkooz1993proper} remains a primary tool for dimensionality reduction, it is restricted by the approximation properties of linear subspaces. For problems with strong nonlinearities, such as turbulent flows, this necessitates a large number of modes to accurately represent the system's dynamics. This limitation has prompted efforts to bridge the gap between linear and nonlinear representations by approximating the system state with a linear projection plus a nonlinear correction term. Notable strategies include quadratic manifolds (QM)~\cite{JAIN201780, rutzmoser2017generalization, barnett2022quadratic, geelen2023operator, SHARMA2023116402, schwerdtner2024greedy, schwerdtner2025online, schwerdtner2025sparseqm, chmiel2025lspg}, which augment linear representations with structured quadratic terms; higher-order polynomial and rational manifolds~\cite{geelen2023learning, geelen2024learning, buchfink2024polymanifolds, KLEIN2025113817, YONG2025117638}, offering more general approximation spaces; and probabilistic manifold decompositions (PMD)~\cite{guo2026nonlinear, guo2026parametricprobabilisticmanifolddecomposition}, which embed residuals into probabilistic latent spaces. Similarly, recent works have employed Gaussian process regression (GPR) and radial basis function (RBF) interpolation~\cite{barnett2023mitigating, ARESDEPARGA2026118443}, deep operator networks (DeepONets)~\cite{codega2026mlclosures}, and reproducing kernel Hilbert spaces (RKHS)~\cite{DIAZ2026118734,diaz2025kernelmanifold} to learn the mapping from the latent state to the correction term. These methods all employ regression to reconstruct the components of the approximation space not captured by the linear basis, offering a data-efficient way to account for the physics discarded during traditional modal truncation while maintaining the original reduced-order dimension. We distinguish these nonlinear representation approaches from classical ROM \emph{closure} modeling. While the latter introduces correction terms to the ROM governing equations to account for truncated modes, the methods discussed here enrich the approximation space itself. For a comprehensive review on modern ROM closures, we refer the reader to~\cite{ahmed2021closures}.

In what follows, we focus on QMs, which offers a compelling balance between nonlinear expressivity and structural transparency. While initial QM implementations use leading POD modes to define the linear subspace~\cite{barnett2022quadratic, geelen2023learning, geelen2024learning}, recent work by Schwerdtner and Peherstorfer has demonstrated that the principal components optimal for linear reconstruction are not necessarily the most effective basis for quadratic representation~\cite{schwerdtner2024greedy,schwerdtner2025online}. In these works, a greedy algorithm is used to select an optimal subset of modes from a larger candidate pool, and it is demonstrated that higher-order modes often contain crucial information for capturing nonlinear interactions. However, such greedy strategies are inherently combinatorial and discrete, restricting the search to fixed singular vector orientations and suffering from factorial complexity as the candidate pool grows. Recent extensions have incorporated dimension-dependent regularization~\cite{ji2026nonlinearmodelorderreduction}, yet the fundamental limitation remains: the subspace is constructed via selection rather than continuous optimization. It should also be noted that the choice of POD coordinates is largely one of practical convenience. Although POD is easy to implement, its optimality is strictly limited to the class of linear approximations, as it minimizes the $\ell_2$ reconstruction error between snapshots and their projection onto a fixed-dimensional subspace. Coordinate transformations have been proposed to optimize the regression landscape~\cite{geelen2023learning, geelen2024learning}, but these formulations remain computationally intensive and lack the scalability required for large-scale applications.

This work proposes a shift from discrete basis selection to continuous subspace identification via Riemannian optimization. Moving beyond the combinatorial selection of POD modes, we frame the search for an optimal basis as an optimization problem on the Stiefel manifold. By searching the space of all orthonormal frames, we identify a subspace specifically well suited for quadratic reconstruction. Our framework restricts the linear subspace to a predefined span of features (the leading POD modes) within which an optimal representation is identified. This allows for the continuous rotation of the basis to minimize quadratic reconstruction error, at a cost independent of the state-space dimension. By exploring arbitrary linear combinations within this feature span, our formulation provides a richer search space than greedy selection and optimizes the synergy between the linear and quadratic terms. We demonstrate superior reconstruction fidelity for a given reduced dimension relative to POD-based QM approaches from~\cite{geelen2024learning} and~\cite{schwerdtner2024greedy} with minimal computational overhead. The \texttt{Pymanopt} library~\cite{townsend2016pymanopt} facilitates Riemannian optimization by leveraging automatic differentiation for all necessary gradient calculations. The primary contributions of this paper include a novel optimization framework for QM construction, the derivation of an objective function whose evaluation does not involve the full state dimension, and a comparison against standard POD and greedily constructed QMs. We demonstrate the robustness of this approach using a large-eddy simulation of a turbulent airfoil wake, showing that optimized subspace rotation consistently yields superior accuracy for a given reduced dimension.

The remainder of this paper is organized as follows. \Cref{sec:dimensionality_reduction} reviews the theoretical foundations of POD and greedy quadratic manifolds. \Cref{sec:riemann} details the proposed Riemannian approach and the derivation of an efficient error metric. \Cref{sec:numerical_experiments} presents numerical results for a turbulent airfoil wake, and \Cref{sec:concluding_remarks} provides concluding remarks.

\section{Dimensionality Reduction}
\label{sec:dimensionality_reduction}

\subsection{Proper Orthogonal Decomposition}
\label{subsec:pod}

Let $\mathbf{s}(t)\in\mathbb{R}^{N}$ be the state of a high-dimensional dynamical process, for example the solution to a spatially discretized PDE model. Following the method of snapshots introduced by Sirovich~\cite{sirovich1987turbulence}, given snapshot observations $\mathbf{s}(t_1),\mathbf{s}(t_2),\ldots,\mathbf{s}(t_K)\in\mathbb{R}^{N}$ of $\mathbf{s}(t)$ at time instances $t_1<t_2<\cdots<t_K$, we form the centered state snapshot matrix
\begin{equation}
    \mathbf{S} = \left[\begin{array}{cccc}
        | & | & & | \\
        (\mathbf{s}(t_1)-\bar{\mathbf{s}}) & (\mathbf{s}(t_2)-\bar{\mathbf{s}}) & \cdots & (\mathbf{s}(t_K)-\bar{\mathbf{s}})\\
        | & | & & |
    \end{array}\right]
    \in\mathbb{R}^{N\times K},
    \label{eq:snapshot_matrix} 
\end{equation}
where $\bar{\mathbf{s}} \in \mathbb{R}^N$ is a fixed reference vector that is usually set to the zero vector, an initial condition, or the mean of the snapshots. Proper orthogonal decomposition (POD) reduces the dimensionality of the state $\mathbf{s}(t)$ by finding an $r$-dimensional subspace $\mathcal{V}_r$ with $r \ll N$ that captures the maximum linear variance of the data. The POD basis for $\mathcal{V}_r$ is given by the $r$ leading left singular vectors $\mathbf{v}_1, \dots, \mathbf{v}_r\in\mathbb{R}^{N}$ of $\mathbf{S}$, corresponding to the largest singular values $\sigma_1 \geq \cdots \geq \sigma_r \geq \dots  \geq \sigma_K \geq 0$ of $\mathbf{S}$.
The state can then be approximated at any time $t$ as a linear combination of the basis vectors plus the reference vector,
\begin{equation}
    \mathbf{s}(t)
    \approx \bar{\mathbf{s}} + \Vr \hat{\mathbf{s}}(t),
    \label{eq:POD_approximation}
\end{equation}
where $\Vr = [~\mathbf{v}_1~~\cdots~~\mathbf{v}_r~]\in\mathbb{R}^{N\times r}$ and $\hat{\mathbf{s}}(t) \in \mathbb{R}^r$ are the state coordinates in the POD basis.
Because $\Vr$ is orthonormal ($\Vr^\top \Vr = \mathbf{I}$), the $L_2$-optimal reduced state is obtained via linear projection,
\begin{equation}
    \hat{\mathbf{s}}(t)
    = \Vr^\top ( \mathbf{s}(t) - \bar{\mathbf{s}}).
    \label{eq:linear_projection}
\end{equation}
The information lost due to truncation is quantified by the relative projection error in the Frobenius norm, which represents the fraction of total energy (variance) not represented in $\mathcal{V}_{r}$:
\begin{equation}
    \frac{\| \mathbf{S} - \Vr \Vr^\top \mathbf{S} \|_F}{\| \mathbf{S} \|_F}
    = \dfrac{\sum_{i=1}^r \sigma_i^2}{\sum_{i=1}^K \sigma_i^2}.
    \label{eq:relative_error}
\end{equation}
The performance of POD-based reduced-order models (ROMs) is inherently tied to the quality of the snapshot data. To ensure predictive accuracy across a desired range of operating conditions, the snapshot matrix $\mathbf{S}$ must be sufficiently representative; in practice, this requires augmenting $\mathbf{S}$ with states from multiple high-fidelity model evaluations, such as those originating from using various initial or boundary conditions. Our methodology applies equally to this situation, but we focus here on single-trajectory data.

\subsection{POD-based Quadratic Manifolds}
\label{subsec:pod_based_qm}

Linear approximation~\cref{eq:POD_approximation} is convenient but can result in significant information loss unless $r$ is relatively large. To address this, we follow~\cite{barnett2022quadratic, geelen2023learning, geelen2024learning} in adding a nonlinear correction term with Kronecker product structure. This yields approximations with a quadratic dependence on the POD coefficients:
\begin{equation}
    \mathbf{s}(t)
    \approx \bar{\mathbf{s}} + \Vr\hat{\mathbf{s}}(t)
    + \hspace{-1em}\underbrace{\Vq \boldsymbol{\Xi}\left( \hat{\mathbf{s}}(t) \otimes \hat{\mathbf{s}}(t) \right)}_{\textstyle \text{\small Quadratic correction term}}\hspace{-1em}, 
\label{eq:quadratic_approximation}
\end{equation}
where $\otimes$ denotes the Kronecker product.
Because $\hat{\mathbf{s}}(t)\otimes\hat{\mathbf{s}}(t)$ contains redundant terms, for computational efficiency (and, later, identifiability) we retain only unique quadratic products so that $\hat{\mathbf{s}}(t)\otimes\hat{\mathbf{s}}(t)$ has dimension $r(r + 1)/2$ (see, e.g.,~\cite{mcquarrie2023thesis}). The columns of the orthonormal matrix $\Vq \in \mathbb{R}^{N \times q}$ are the next available higher-order left singular vectors of the snapshot data, $\mathbf{v}_{r+1},\mathbf{v}_{r+2},\ldots,\mathbf{v}_{r+q}$. Finally, $\boldsymbol{\Xi} \in \mathbb{R}^{q \times r(r+1)/2}$ denotes a weighting matrix that must be inferred from available data.
The approach in~\cite{geelen2024learning} determines $\boldsymbol{\Xi}$ through a linear least-squares problem,
\begin{equation}
    \min_{\boldsymbol{\Xi}} \left\{\left\|
        (\mathbf{I} - \Vr \Vr^\top ) \mathbf{S} - \Vq \boldsymbol{\Xi}\left( (\Vr^\top \mathbf{S}) \odot (\Vr^\top \mathbf{S}) \right)
    \right\|_F^2 + \gamma \| \boldsymbol{\Xi} \|_F^2\right\},
    \label{eq:qm_least_squares}
\end{equation}
where $\mathbf{I}$ is the $N\times N$ identity matrix, $\odot$ denotes the Khatri–Rao (columnwise Kronecker) product, and $\gamma \geq 0$ is a scalar regularization parameter.

Several recent studies on quadratic manifold approximation~\cite{DIAZ2026118734, ji2026nonlinearmodelorderreduction} treat the product $\Vq \boldsymbol{\Xi}$ as a single matrix, but the factorized structure from \cref{eq:quadratic_approximation} offers several distinct advantages. First, by the virtue of the properties of singular vectors, the columns of $\Vr$ are orthogonal to those of $\Vq$ (that is, $\Vr^\top \Vq = \mathbf{0}$). This means that least-squares problem~\cref{eq:qm_least_squares} has only $q \times r(r+1)/2$ unknowns rather than $N\times r(r+1)/2$, paving the way to truly large-scale applications~\cite{geelen2024learning}. Second, it significantly streamlines the implementation of a greedy basis vector selection algorithm, to be discussed next.

\subsection{Greedily Constructed Quadratic Manifolds}
\label{subsec:greedy_qm}

The greedy selection method from Schwerdtner and Peherstorfer~\cite{schwerdtner2024greedy} offers a robust alternative for constructing the quadratic manifold approximation space.
Instead of selecting the \emph{principal} left singular vectors $\mathbf{v}_1, \dots, \mathbf{v}_{r+q}$ of the centered state snapshot matrix $\mathbf{S}$, the procedure sets a fixed number of candidates $m$ with $r \ll m \le K$ and iteratively selects a \emph{non-sequential} set of indices $j_1, \dots, j_r\in\{1,\ldots,m\}$ such that $\{\mathbf{v}_{j_1},\ldots,\mathbf{v}_{j_r}\}$ and $\{\mathbf{v}_{j_{r+1}},\ldots,\mathbf{v}_{j_m}\}$ are optimal in a greedy sense for the linear and quadratic parts, respectively, of approximation \Cref{eq:quadratic_approximation}.
By selecting optimally from the larger pool of $m \gg r$ left singular vectors of $\mathbf{S}$, this strategy effectively captures the necessary physics by allowing the inclusion of high-order modes that traditional linear energy criteria might otherwise discard.

We detail the procedure, adapted to match the approximation form \cref{eq:quadratic_approximation}, as follows.
Let $i = 1,\dots, r$ be the iteration counter variable of the greedy selection and $\mathbf{V}_{\!(i)} = [~\mathbf{v}_{j_1}~~\mathbf{v}_{j_2}~~\cdots~~\mathbf{v}_{j_i}~] \in \mathbb{R}^{N \times i}$ be the basis matrix associated with the linear part of \cref{eq:quadratic_approximation} at iteration $i$. Additionally, let $\mathbf{V}_{\!(!i)} = [~\mathbf{v}_{j_{i+1}}~~\cdots~~\mathbf{v}_{j_{m}}~]  \in \mathbb{R}^{N \times (m-i)}$ be the basis matrix associated with the quadratic correction in \cref{eq:quadratic_approximation}, comprised of the modes \emph{not} selected in $\mathbf{V}_{\!(i)}$ after $i$ iterations. Note that the matrix $[~\mathbf{V}_{\!(i)}~~\mathbf{V}_{\!(!i)}~] \in \mathbb{R}^{N \times m}$ contains all $m$ candidate vectors $\mathbf{v}_1,\ldots,\mathbf{v}_m$. At iteration $i$, we select the left singular vector $\mathbf{v}_{j_i}$ with index $j_i$ that minimizes the representation error:
\begin{equation}
    \min_{j_i\in\{1,\dots,m\}} \min_{\boldsymbol{\Xi}} \left\{
        \left\|
            \left(\mathbf{I} - \mathbf{V}_{\!(i)} \mathbf{V}_{\!(i)}^\top \right)\mathbf{S}
            - \mathbf{V}_{\!(!i)} \boldsymbol{\Xi}\left( (\mathbf{V}_{\!(i)}^\top \mathbf{S}) \odot (\mathbf{V}_{\!(i)}^\top \mathbf{S}) \right)
        \right\|_F^2
        + \gamma \| \boldsymbol{\Xi} \|_F^2 \right\},
    \label{eq:greedy}
\end{equation}
with the basis matrix being updated at every iteration $i$ as
\begin{equation}
    \mathbf{V}_{\!(i)}
    = [~\mathbf{V}_{\!(i-1)}~~\mathbf{v}_{j_i}~] \in \mathbb{R}^{N \times i}.
    \label{eq:greedy2}
\end{equation}
The algorithm starts at iteration $i = 0$ with the subspace $\mathbf{V}_{\!(0)}$ that contains only the zero vector. After $r$ iterations, we obtain the basis matrices $\Vr := \mathbf{V}_{\!(r)} \in \mathbb{R}^{N \times r}$ and $\Vq := \mathbf{V}_{\!(!r)} \in \mathbb{R}^{N \times q}$ in \cref{eq:quadratic_approximation}, where $q = m-r$, and compute the corresponding weight matrix $\boldsymbol{\Xi}$ by solving \cref{eq:qm_least_squares} to complete the approximation.

In contrast to \cite{schwerdtner2024greedy}, the computational burden of our Greedy formulation can (and should) be further reduced by restricting the nonlinear expansion to only $r+q < m $ modes. Such a truncation, supported by~\cite{geelen2024learning} as well as our results in \Cref{sec:numerical_experiments}, strategically filters modes that are not amenable to quadratic approximation. The parametrization in \cref{eq:quadratic_approximation} facilitates such a truncation. This modification is key to preserving a parsimonious model structure and guaranteeing that QM approaches remain computationally tractable in many problem settings.
 
\section{A Riemannian Optimization Perspective}
\label{sec:riemann}

\subsection{Towards Optimal Quadratic Manifolds}

Nonlinear approximations based on left singular vectors alone can miss information that is essential for the quadratic correction term in \cref{eq:quadratic_approximation} to be effective. Rather than using a combinatorial selection of POD modes, we structure the search for an optimal basis as an optimization problem over the space of all orthonormal frames with the goal of identifying a subspace that is optimally suited for quadratic approximation. We focus on the following constrained optimization problem, formulated as a regularized least-squares objective:
\begin{equation}
\begin{aligned}
    & \min_{\boldsymbol{\Xi}, \Vr, \Vq} \left\{ \left\| (\mathbf{I} - \Vr \Vr^\top ) \mathbf{S} - \Vq \boldsymbol{\Xi}\left( (\Vr^\top \mathbf{S}) \odot (\Vr^\top \mathbf{S}) \right)
    \right\|_F^2 + \gamma \| \boldsymbol{\Xi} \|_F^2 \right\} \\
    & \text{subject to} \quad [~\Vr~~\Vq~] \in \text{St}(N, r+q),
    \label{eq:qm_least_squares2}
\end{aligned}
\end{equation}
where the \emph{Stiefel} manifold $\text{St}(N, r+q)$ is set of $N \times (r+q)$ matrices with orthonormal columns. Because the design space scales with the dimensionality $N$ of the state space, numerically solving \cref{eq:qm_least_squares2} is generally computationally intractable. To overcome this challenge, we restrict the search space to a predefined feature span comprising the leading $m \ll N$ left singular vectors of the data. This is achieved through the factorizations
\begin{equation}
    \Vr = \tilde{\mathbf{V}} \mathbf{Q}_r,
    \qquad
    \Vq = \tilde{\mathbf{V}} \mathbf{Q}_q,
\end{equation}
where the matrix $\tilde{\mathbf{V}} = [\tilde{\mathbf{v}}_1, \dots, \tilde{\mathbf{v}}_m ] \in \mathbb{R}^{N \times m}$ assembles $m$ left singular vectors. We constrain the matrices $\mathbf{Q}_r \in \mathbb{R}^{m \times r}$ and  $\mathbf{Q}_q \in \mathbb{R}^{m \times q}$ to be orthonormal with mutually orthogonal column spaces, that is,
\begin{equation}
    \mathbf{Q}_r^\top\mathbf{Q}_r = \mathbf{I}_r,
    \qquad
    \mathbf{Q}_r^\top\mathbf{Q}_q = \mathbf{0},
    \qquad
    \mathbf{Q}_q^\top\mathbf{Q}_q = \mathbf{I}_q.
    \label{eq:q-orthonormality}
\end{equation}
Consequently, every $j$th column vector in $\Vr$ and $\Vq$, represented by $\mathbf{v}_{r, j}$ and $\mathbf{v}_{q, j}$, can be written as a linear combination of candidate basis vectors:
\begin{equation}
\begin{aligned} 
    \mathbf{v}_{r,j}
    &= Q^{(1,j)}_r\tilde{\mathbf{v}}_1
    + Q^{(2,j)}_r\tilde{\mathbf{v}}_2
    + \dots + Q^{(m,j)}_r\tilde{\mathbf{v}}_m,
    \\
    \mathbf{v}_{q,j}
    &= Q^{(1,j)}_q\tilde{\mathbf{v}}_1
    + Q^{(2,j)}_q\tilde{\mathbf{v}}_2
    + \dots + Q^{(m,j)}_q \tilde{\mathbf{v}}_m,
\end{aligned}
\end{equation}
where $Q^{(i,j)}_r$ and $Q^{(i,j)}_q$ denote the $(i,j)$ elements of the rotation matrices $\mathbf{Q}_r$ and $\mathbf{Q}_q$, respectively. The vectors $\mathbf{v}_{r,1},\ldots\mathbf{v}_{r,r}\in\mathbb{R}^{N}$ and $ \mathbf{v}_{q,1},\ldots,\mathbf{v}_{q,q}\in\mathbb{R}^{N}$ will be referred to as \emph{optimized} basis vectors (also  modes). Effectively, each of the $r$ and $q$ optimized vectors is reconstructed from the $m$-dimensional candidate span. Importantly, we must maintain
\begin{equation}
     r + q \leq m
     \label{eq:linear_independence}
\end{equation}
to ensure the optimized basis vectors remain linearly independent since each $\mathbf{v}_{r,j}$ and each $\mathbf{v}_{q,j}$ is contained in the $m$-dimensional column span of $\tilde{\mathbf{V}}$.

With these definitions in hand, the columns of $\Vr$ and $\Vq$ form an orthonormal set:
\begin{equation}
    \begin{aligned}
    \Vr^\top \Vr = (\mathbf{Q}_r^\top \tilde{\mathbf{V}}^\top) (\tilde{\mathbf{V}} \mathbf{Q}_r) = \mathbf{Q}_r^\top \mathbf{I}_{m} \mathbf{Q}_r = \mathbf{I}_r, \\
    \Vq^\top \Vq = (\mathbf{Q}_q^\top \tilde{\mathbf{V}}^\top) (\tilde{\mathbf{V}}\mathbf{Q}_q) = \mathbf{Q}_q^\top \mathbf{I}_{m} \mathbf{Q}_q = \mathbf{I}_q.
    \end{aligned}
\end{equation}
Moreover, the columns of $\Vr$ are strictly orthogonal to the columns of $\Vq$:
\begin{equation}
    \Vr^\top \Vq = (\mathbf{Q}_r^\top \tilde{\mathbf{V}}^\top) (\tilde{\mathbf{V}}\mathbf{Q}_q) = \mathbf{Q}_r^\top \mathbf{I}_m \mathbf{Q}_q = \mathbf{Q}_r^\top \mathbf{Q}_q = \mathbf{0}.
\end{equation}
By using a precomputable, fixed high-dimensional basis, we have shifted the optimization burden from the high-dimensional state space to the determination of linear coefficients $Q_r^{(i,j)}, Q_q^{(i,j)}$. This approach ensures that the high-dimensional matrices $\mathbf{V}_r$ and $\mathbf{V}_q$ do not need explicit updating during the numerical search. The optimization problem can now be written as
\begin{equation}
    \begin{aligned}
    &\min_{\boldsymbol{\Xi}, \mathbf{Q}_r, \mathbf{Q}_q} \mathcal{J}(\boldsymbol{\Xi}, \mathbf{Q}_r, \mathbf{Q}_q) \\ & \text{subject to} \quad [~\mathbf{Q}_r~~\mathbf{Q}_q~] \in \text{St}(m, r+q),
    \label{eq:qm_least_squares3}
\end{aligned}
\end{equation}
where
\begin{equation}
    \mathcal{J}(\boldsymbol{\Xi}, \mathbf{Q}_r, \mathbf{Q}_q) = \left\|
            (\mathbf{I} - \tilde{\mathbf{V}} \mathbf{Q}_r \mathbf{Q}_r^\top  \tilde{\mathbf{V}}^\top ) \mathbf{S}
            - \tilde{\mathbf{V}} \mathbf{Q}_q \boldsymbol{\Xi}  \mathbf{W}
        \right\|_F^2
        + \gamma \| \boldsymbol{\Xi} \|_F^2,
        \label{eq:objective}
\end{equation}
with
\begin{equation}
     \mathbf{W} = \mathbf{W}(\mathbf{Q}_r) = (\mathbf{Q}_r^\top  \tilde{\mathbf{V}}^\top \mathbf{S}) \odot (\mathbf{Q}_r^\top  \tilde{\mathbf{V}}^\top \mathbf{S}) \in \mathbb{R}^{r(r+1)/2 \times K}.
\end{equation}
We can then also exploit the fact that the optimal coefficients $\boldsymbol{\Xi} = \boldsymbol{\Xi}(\mathbf{Q}_r,\mathbf{Q}_q)$ are dependent only on $\mathbf{Q}_r$ and $\mathbf{Q}_q$ through the solution to the linear least-squares problem \cref{eq:qm_least_squares}, which satisfies
\begin{equation}
    \boldsymbol{\Xi} \left( \mathbf{W} \mathbf{W}^\top + \gamma \mathbf{I} \right) = \mathbf{Q}_q^\top \tilde{\mathbf{V}}^\top \mathbf{S}  \mathbf{W}^\top.
    \label{eq:least_squares}
\end{equation}
By substituting the solution $\boldsymbol{\Xi}$ back into the objective function, we reduce the optimization task to a problem solely dependent on $\mathbf{Q}_r$ and $\mathbf{Q}_q$:
\begin{equation}
    \begin{aligned}
    &\min_{\mathbf{Q}_r, \mathbf{Q}_q} \mathcal{J}( \mathbf{Q}_r, \mathbf{Q}_q) \\ & \text{subject to} \quad [~\mathbf{Q}_r~~\mathbf{Q}_q~] \in \text{St}(m, r+q).
    \label{eq:qm_least_squares4}
\end{aligned}
\end{equation}
It is of the utmost importance that the cost of numerically evaluating the objective function remains independent of the state-space dimension $N$. To ensure computational efficiency, we reformulate the objective to eliminate any operations that scale with the full system size, as detailed in the following section.

\subsection{Reformulating the Optimization Problem}

Through algebraic manipulation of the objective in \cref{eq:objective}, we identify the individual terms that may be efficiently computed offline. This significantly reduces the overhead required for both function evaluations and gradient computations. We start with the compact form
\begin{equation}
    \left\| \mathbf{S} - \tilde{\mathbf{V}} \mathbf{Q}_r \mathbf{Q}_r^\top  \tilde{\mathbf{V}}^\top \mathbf{S} - \tilde{\mathbf{V}} \mathbf{Q}_q \boldsymbol{\Xi} \mathbf{W} \right\|_F^2 + \gamma \| \boldsymbol{\Xi} \|_F^2.
\end{equation}
Utilizing the trace representation of the Frobenius norm, we obtain
\begin{equation}
\begin{aligned}
& \| \mathbf{S} - \tilde{\mathbf{V}}\mathbf{Q}_r \mathbf{Q}_r^\top \tilde{\mathbf{V}}^\top \mathbf{S}
- \tilde{\mathbf{V}}\mathbf{Q}_q \boldsymbol{\Xi} \mathbf{W} \|_F^2 + \gamma \| \boldsymbol{\Xi} \|_F^2
\\
& \qquad = \text{Tr}\left(
(\mathbf{S} - \tilde{\mathbf{V}}\mathbf{Q}_r \mathbf{Q}_r^\top \tilde{\mathbf{V}}^\top \mathbf{S}
- \tilde{\mathbf{V}}\mathbf{Q}_q \boldsymbol{\Xi} \mathbf{W})^\top (\mathbf{S} - \tilde{\mathbf{V}}\mathbf{Q}_r \mathbf{Q}_r^\top \tilde{\mathbf{V}}^\top \mathbf{S}
- \tilde{\mathbf{V}}\mathbf{Q}_q \boldsymbol{\Xi} \mathbf{W})
\right) + \gamma \text{Tr} \left( \boldsymbol{\Xi}^{\top} \boldsymbol{\Xi} \right) \\
& \qquad = \text{Tr}\Big(
\mathbf{S}^\top\mathbf{S}
- \mathbf{S}^\top \tilde{\mathbf{V}}\mathbf{Q}_r \mathbf{Q}_r^\top \tilde{\mathbf{V}}^\top \mathbf{S}
- \mathbf{S}^\top \tilde{\mathbf{V}}\mathbf{Q}_q \boldsymbol{\Xi} \mathbf{W} - \mathbf{S}^\top \tilde{\mathbf{V}}\mathbf{Q}_r \mathbf{Q}_r^\top \tilde{\mathbf{V}}^\top \mathbf{S}
+ \mathbf{S}^\top \tilde{\mathbf{V}}\mathbf{Q}_r \mathbf{Q}_r^\top \tilde{\mathbf{V}}^\top
  \tilde{\mathbf{V}}\mathbf{Q}_r \mathbf{Q}_r^\top \tilde{\mathbf{V}}^\top \mathbf{S} \\
&\qquad \qquad
+ \mathbf{S}^\top \tilde{\mathbf{V}}\mathbf{Q}_r \mathbf{Q}_r^\top \tilde{\mathbf{V}}^\top
  \tilde{\mathbf{V}}\mathbf{Q}_q \boldsymbol{\Xi} \mathbf{W}
+ \mathbf{W}^\top \boldsymbol{\Xi}^\top \mathbf{Q}_q^\top \tilde{\mathbf{V}}^\top
  \tilde{\mathbf{V}}\mathbf{Q}_r \mathbf{Q}_r^\top \tilde{\mathbf{V}}^\top \mathbf{S} - \mathbf{W}^\top \boldsymbol{\Xi}^\top \mathbf{Q}_q^\top \tilde{\mathbf{V}}^\top \mathbf{S}
\\
& \qquad \qquad + \mathbf{W}^\top \boldsymbol{\Xi}^\top \mathbf{Q}_q^\top \tilde{\mathbf{V}}^\top
  \tilde{\mathbf{V}}\mathbf{Q}_q \boldsymbol{\Xi} \mathbf{W} + \gamma \boldsymbol{\Xi}^{\top} \boldsymbol{\Xi} \Big).
\end{aligned}
\label{eq25}
\end{equation}
Using the orthogonality conditions~\cref{eq:q-orthonormality} and recalling $\tilde{\mathbf{V}}^\top\tilde{\mathbf{V}}= \mathbf{I}$, the above expression simplifies to 
\begin{equation}
\begin{aligned}
    & \| \mathbf{S} - \tilde{\mathbf{V}}\mathbf{Q}_r \mathbf{Q}_r^\top \tilde{\mathbf{V}}^\top  \mathbf{S} - \tilde{\mathbf{V}}\mathbf{Q}_q \boldsymbol{\Xi} \mathbf{W} \|_F^2
    + \gamma \| \boldsymbol{\Xi} \|_F^2
    \\
    & \qquad = \text{Tr} \left(
        \mathbf{S}^\top\mathbf{S}
        - \mathbf{S}^\top \tilde{\mathbf{V}}\mathbf{Q}_r \mathbf{Q}_r^\top \tilde{\mathbf{V}}^\top  \mathbf{S}
        - \mathbf{S}^\top \tilde{\mathbf{V}}\mathbf{Q}_q \boldsymbol{\Xi} \mathbf{W}
        - \mathbf{W}^\top \boldsymbol{\Xi}^\top \mathbf{Q}_q^\top \tilde{\mathbf{V}}^\top \mathbf{S}
        + \mathbf{W}^\top \boldsymbol{\Xi}^\top \boldsymbol{\Xi} \mathbf{W}
        + \gamma \boldsymbol{\Xi}^{\top} \boldsymbol{\Xi} \right)
    \\
    & \qquad = \underbrace{\| \mathbf{S} \|_F^2}_\text{precomputable} 
    \hspace{-1em} - \| \mathbf{Q}_r^\top (\hspace{-1em}\underbrace{\tilde{\mathbf{V}}^\top \mathbf{S}}_\text{precomputable} \hspace{-1em}) \|_F^2
    + \| \boldsymbol{\Xi}\mathbf{W} \|_F^2 + \gamma \| \boldsymbol{\Xi} \|_F^2
    - 2 \text{Tr} ( ( \hspace{-1em} \underbrace{\tilde{\mathbf{V}}^\top \mathbf{S}}_\text{precomputable}\hspace{-1em})^\top \mathbf{Q}_q \boldsymbol{\Xi} \mathbf{W}).
\end{aligned}
\label{eq25-1}
\end{equation}
With $\tilde{\mathbf{V}}^\top \mathbf{S}$ and $\| \mathbf{S} \|_F^2$ precomputed, the full objective function can be evaluated without reference to the high dimension $N$. The problem now resides solely in the lower-dimensional feature space of the predefined candidate vectors, where we recall that $r + q \leq m \leq K \ll N$. Because $\| \mathbf{S} \|_F^2$ does not depend on $\mathbf{Q}_r$ or $\mathbf{Q}_q$, it may be omitted from the minimization problem. If we then introduce 
\begin{equation}
\tilde{\mathbf{S}} = \tilde{\mathbf{V}}^\top \mathbf{S} \in \mathbb{R}^{m \times K}
\end{equation}
as the representation of snapshot data $\mathbf{S}$ in the candidate basis $\tilde{\mathbf{V}}$, optimization problem \cref{eq:qm_least_squares4} can be written as
\begin{equation}
\begin{aligned}
    &\min_{\mathbf{Q}_r, \mathbf{Q}_q} \left\{
        - \| \mathbf{Q}_r^\top \tilde{\mathbf{S}} \|_F^2
        + \| \boldsymbol{\Xi}\mathbf{W} \|_F^2
        + \gamma \| \boldsymbol{\Xi} \|_F^2
        - 2 \text{Tr} ( \tilde{\mathbf{S}}^\top \mathbf{Q}_q \boldsymbol{\Xi} \mathbf{W})
    \right\} \\ & \text{subject to} \quad [~\mathbf{Q}_r~~\mathbf{Q}_q~] \in \text{St}(m, r+q).
    \label{eq:qm_least_squares5}
\end{aligned}
\end{equation}

\subsection{Method Overview} 

\begin{algorithm}[!tbp]
\caption{\texttt{FastQM}: Construct quadratic manifolds through Riemannian optimization}
\label{alg:fastqm_algorithm}
\begin{algorithmic}[1]
    \Require Centered snapshot matrix $\mathbf{S} \in \mathbb{R}^{N \times K}$;  dimensionalities $r,q$ in \cref{eq:quadratic_approximation}; number of candidate modes $m$
    \State $\mathbf{X}\boldsymbol{\Sigma}\mathbf{Z}^\top \gets \operatorname{SVD}(\mathbf{S})$
        \Comment{Compute thin SVD of centered snapshot data.}
    \State $\tilde{\mathbf{V}} \gets \mathbf{X}_{:,:m}$
        \Comment{Select $m$ candidate modes (leading left singular vectors).}
    \State $\tilde{\mathbf{S}} \gets \tilde{\mathbf{V}}^\top\mathbf{S} \in \mathbb{R}^{m \times K}$
        \Comment{Project centered snapshots to column span of $\tilde{\mathbf{V}}$.}
    \State Solve the following orthogonality-constrained optimization problem using \texttt{Pymanopt}~\cite{townsend2016pymanopt}:
    \begin{equation*}
    \min_{\mathbf{Q}_r, \mathbf{Q}_q}  \mathcal{J}(\mathbf{Q}_r, \mathbf{Q}_q); \quad \text{subject to} \quad [~\mathbf{Q}_r~~\mathbf{Q}_q~] \in \text{St}(m, r+q),
    \end{equation*}
    where
    \begin{equation*}
    \mathcal{J}(\mathbf{Q}_r, \mathbf{Q}_q) = - \| \mathbf{Q}_r^\top \tilde{\mathbf{S}} \|_F^2 + \| \boldsymbol{\Xi}\mathbf{W} \|_F^2 + \gamma \| \boldsymbol{\Xi} \|_F^2 - 2 \text{Tr} ( \tilde{\mathbf{S}}^\top \mathbf{Q}_q \boldsymbol{\Xi} \mathbf{W})
    \end{equation*}
    \item[]and where $\boldsymbol{\Xi} = \boldsymbol{\Xi}(\mathbf{Q}_r,\mathbf{Q}_q)$ solves the linear least-squares problem \cref{eq:least_squares}.
    \Ensure Rotation matrices $\mathbf{Q}_r$ and $\mathbf{Q}_q$; coefficient matrix $\boldsymbol{\Xi}(\mathbf{Q}_r,\mathbf{Q}_q)$
\end{algorithmic}
\end{algorithm}

\begin{figure}[!tbp]
\centering
\begin{subfigure}[t]{.9\linewidth}
    \centering
    \includegraphics[width=\linewidth]{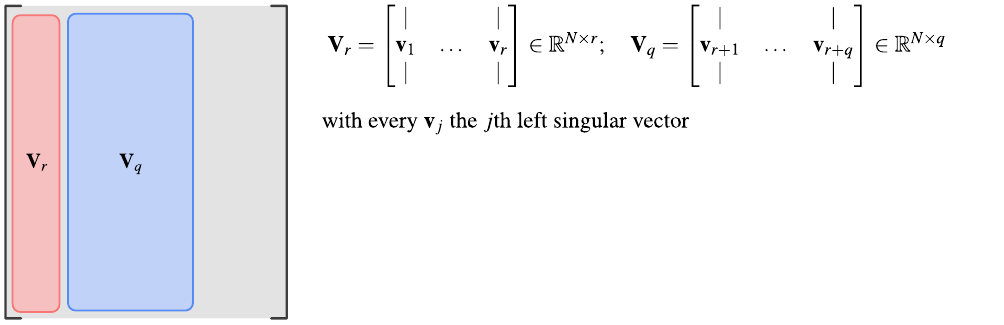}
    \caption{POD-based QM strategy (see~\Cref{subsec:pod_based_qm})}
\end{subfigure}\\[1em]
\begin{subfigure}[t]{.9\linewidth}
    \centering
    \includegraphics[width=\linewidth]{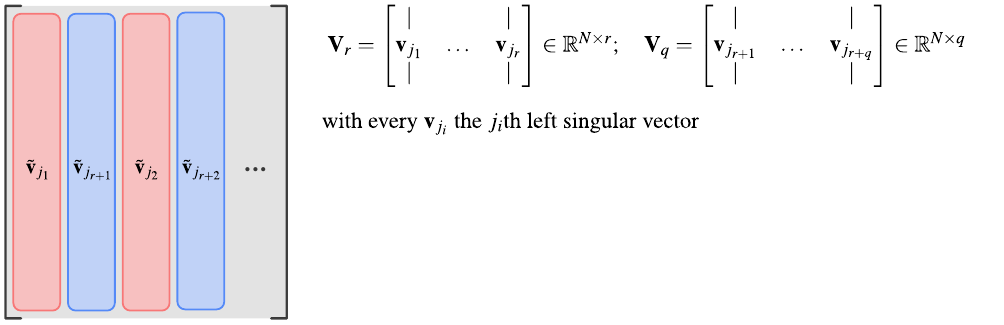}
    \caption{Greedy-based QM strategy  (see~\Cref{subsec:greedy_qm})}
\end{subfigure}\\[1em]
\begin{subfigure}[t]{.9\linewidth}
    \centering
    \includegraphics[width=\linewidth]{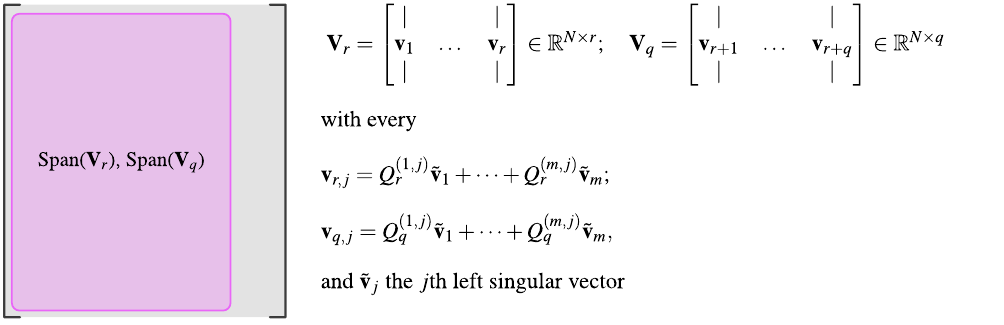}
    \caption{Proposed Riemannian approach}
\end{subfigure}
\caption{Building quadratic manifold approximations of the form \cref{eq:quadratic_approximation} using (a) leading left singular vectors~\cite{geelen2023learning, geelen2024learning}, (b) Greedily selected left singular vectors~\cite{schwerdtner2024greedy} or (c) the proposed Riemannian approach which amounts to finding two spanning sets contained within the column space spanned by $m$ candidate left singular vectors.}
\label{fig:overview}
\end{figure}

By shifting from discrete mode selection to continuous subspace identification via Riemannian optimization, the \texttt{FastQM} algorithm efficiently learns quadratic manifolds by structuring the search for an optimal basis on the Stiefel manifold rather than performing a combinatorial search across a fixed set of singular vectors. To maintain numerical feasibility, the search space is restricted to a predefined feature span of the leading $m$ left singular vectors of the data matrix. Within this span, we seek orthogonal matrices that identify a subspace that is optimal for the purpose of learning quadratic manifold approximations. Our implementation is fundamentally anchored by \texttt{Pymanopt}~\cite{townsend2016pymanopt}, which serves as a critical enabler by leveraging automatic differentiation to evaluate complex manifold gradients with machine precision; by automating gradient computations, we eliminate the potential for human error inherent in manual derivations. Furthermore, this automated approach ensures that the optimization remains numerically stable even as dimension of the search space increases for problems with high complexity. For the numerical backend, we utilize \texttt{JAX}~\cite{jax2018github} to provide hardware-accelerated linear algebra, ensuring the optimization remains efficient as the candidate pool size increases even for large-scale problem settings. The overall process of the proposed quadratic manifold construction using Riemannian optimization is described in \Cref{alg:fastqm_algorithm}. The proposed subspace identification strategy is contrasted with traditional POD-based and greedy  QM constructions in \Cref{fig:overview}.

\subsection{Motivational Example}

To demonstrate why POD modes are often restrictive in quadratic approximation and to offer deeper insight into the proposed Riemannian approach, we consider the two-dimensional trajectory $\mathbf{s}(t) = [~s_1(t)~~s_2(t)~]^\top =  [~t~~t^2~]^\top$ over the time interval $t \in [-0.75, 1.25]$. The setup for this problem is adapted from \cite{KORONAKI2024112910}. Despite the low dimensionality, we attempt to approximate the trajectory using only a single degree of freedom $\hat{s}(t)$. We define the snapshot matrix $\mathbf{S} \in \mathbb{R}^{2 \times 25}$ using 25 uniformly sampled $t$-values. The reference state $\bar{\mathbf{s}}$ is chosen as $\mathbf{s}(0) = \mathbf{0}$.

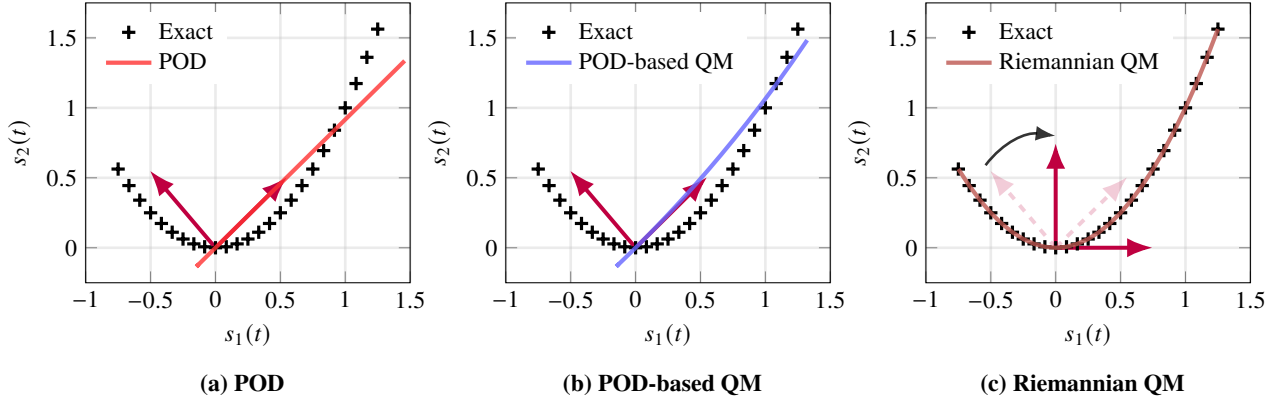
\begin{figure}[!tbp]
\vspace{1em}
    \centering \small
    \hspace{2em}
    \begin{subfigure}[t]{0.26\linewidth}
        \hspace{-3.5em}
        \begin{tikzpicture}
        \begin{axis}[
            set layers,
            width=\linewidth,
            scale only axis,
            xlabel={$s_1(t)$},
            ylabel={$s_2(t)$},
            xmin=-1, xmax=1.5,
            ymin=-0.25, ymax=1.75,
            mark size=1pt,
            grid=both,   
            major grid style={line width=1pt, opacity=0.3},
            xtick={-1, -0.5, ..., 1.5}, 
            ytick={},
            legend pos=north west,
            legend style={
                draw=none,
                legend cell align=left
            }            
            ]
            
        \addplot[very thick, 
            only marks,
            mark=+, 
            draw=black,
            mark options={fill=black},
            mark size=2.5pt,
            on layer=axis foreground
        ] table [
            header=false,
            x index=0,
            y index=1,
            col sep=comma
        ]{Data/2D_Example.csv};
        \addlegendentry{Exact}

        \addplot[
            ultra thick,
            mark=none,
            opacity=0.8,
            draw=red!80,
            mark size=2.5pt,
            on layer=axis foreground
        ] table [
            header=false,
            x index=2,
            y index=3,
            col sep=comma
        ]{Data/2D_Example.csv};
        \addlegendentry{POD}

        \draw[-{Latex[purple]}, purple, line width=1.5] (0,0) -- (0.7372444526469064*0.75,0.6756260926291728*0.75);
        \draw[-{Latex[purple]}, purple, line width=1.5] (0,0) -- (-0.6756260926291728*0.75,0.7372444526469064*0.75);

        \end{axis}
        \end{tikzpicture}
        \caption{POD}
        \label{fig:2D_POD}
    \end{subfigure}
    \hspace{3.5em}
    \begin{subfigure}[t]{0.26\linewidth}
        \hspace{-3.5em}
        \begin{tikzpicture}
        \begin{axis}[
            set layers,
            width=\linewidth,
            scale only axis,
            xlabel={$s_1(t)$},
            ylabel={$s_2(t)$},
            xmin=-1, xmax=1.5,
            ymin=-0.25, ymax=1.75,
            mark size=1pt,
            grid=both,   
            major grid style={line width=1pt, opacity=0.3},
            xtick={-1, -0.5, ..., 1.5}, 
            ytick={},
            legend pos=north west,
            legend style={
                draw=none,
                legend cell align=left
            }            
            ]
            
        \addplot[very thick, 
            only marks,
            mark=+, 
            draw=black,
            mark options={fill=black},
            mark size=2.5pt
        ] table [
            header=false,
            x index=0,
            y index=1,
            col sep=comma
        ]{Data/2D_Example.csv};
        \addlegendentry{Exact}
    
        \addplot[
            ultra thick,
            mark=none,
            opacity=0.8,
            draw=blue!60,
            mark size=2.5pt,
            on layer=axis foreground
        ] table [
            header=false,
            x index=4,
            y index=5,
            col sep=comma
        ]{Data/2D_Example.csv};
        \addlegendentry{POD-based QM}

        \draw[-{Latex[purple]}, purple, line width=1.5] (0,0) -- (0.7372444526469064*0.75,0.6756260926291728*0.75);
        \draw[-{Latex[purple]}, purple, line width=1.5] (0,0) -- (-0.6756260926291728*0.75,0.7372444526469064*0.75);
        
        \end{axis}
        \end{tikzpicture}
        \caption{POD-based QM}
        \label{fig:2D_QM}
    \end{subfigure}
    \hspace{3.5em}
    \begin{subfigure}[t]{0.26\linewidth}
        \hspace{-3.5em}
        \begin{tikzpicture}
        \begin{axis}[
            set layers, 
            width=\linewidth,
            scale only axis,
            xlabel={$s_1(t)$},
            ylabel={$s_2(t)$},
            xmin=-1, xmax=1.5,
            ymin=-0.25, ymax=1.75,
            mark size=1pt,
            legend pos=north east,
            grid=both,   
            major grid style={line width=1pt, opacity=0.3},
            xtick={-1, -0.5, ..., 1.5}, 
            ytick={},
            legend pos=north west,
            legend style={
                draw=none,
                legend cell align=left
            }            
            ]
            
        \addplot[very thick, 
            only marks,
            mark=+, 
            draw=black,
            mark options={fill=black},
            mark size=2.5pt
        ] table [
            header=false,
            x index=0,
            y index=1,
            col sep=comma
        ]{Data/2D_Example.csv};
        \addlegendentry{Exact}
    
        \addplot[
            ultra thick,
            mark=none,
            opacity=0.8,
            draw=Maroon!80,
            mark size=2.5pt,
            on layer=axis foreground
        ] table [
            header=false,
            x index=6,
            y index=7,
            col sep=comma
        ]{Data/2D_Example.csv};
        \addlegendentry{Riemannian QM}

        \draw[-{Latex[purple]}, purple, line width=1.5, dashed, opacity=0.2] (0,0) -- (0.7372444526469064*0.75,0.6756260926291728*0.75);
        \draw[-{Latex[purple]}, purple, line width=1.5, dashed, opacity=0.2] (0,0) -- (-0.6756260926291728*0.75,0.7372444526469064*0.75);

        \draw [-{Latex}, black, line width=1, opacity=0.8] (-0.6756260926291728*0.8,0.7372444526469064*0.8) to [bend left] (-0.0003975556949919269*0.8,0.9999999209747317*0.8);
        
        \draw[-{Latex[purple]}, purple, line width=1.5] (0,0) -- (0.9999999209747317*0.75,0.0003975556949919269*0.75);
        \draw[-{Latex[purple]}, purple, line width=1.5] (0,0) -- (-0.0003975556949919269*0.75,0.9999999209747317*0.75);

        \end{axis}
        \end{tikzpicture}
        \caption{Riemannian QM}
        \label{fig:2D_R}
    \end{subfigure}
    \caption{Comparison of the approximations produced by (a) POD; (b) POD-based QM approximations; (c) and QM approximations obtained via Riemannian optimization.}
\end{figure}

\Cref{fig:2D_POD} represents the POD approximation with the red curve, with the left singular vectors indicated by the two red arrows. This approximation fails to capture the intrinsic curvature of the data, resulting in a substantial relative error of 37.32\%. The POD-based QM approach yields an approximation that is no longer strictly tied to the direction of the first singular vector, enabling it to capture solution features outside of that single dimension (see \Cref{fig:2D_QM}). However, while the relative error improves slightly to 36.59\%, it is clear that the left singular vectors provide an inadequate basis for QM approximation.

It is evident that by allowing the basis to \emph{rotate}, the vectors may positioned such that the linear and quadratic terms integrate more effectively, better capturing the curvature of the solution manifold. Indeed, an ideal basis spanned by $[~1~~0~]^\top$ and $[~0~~1~]^\top$ would allow for the exact reconstruction of every point along the trajectory $\mathbf{s}(t)$. Applying the proposed Riemannian approach, the required rotation for the original POD basis is determined by solving for $\mathbf{Q}_r$ and $\mathbf{Q}_q$. Solving the optimization problem yields the approximation
\begin{equation}
    \mathbf{s}(t)
    \approx \underbrace{\left[\begin{array}{rr}
    0.737 & -0.676\\ 
    0.676 &  0.737
    \end{array}\right]}_{\textstyle \text{\small$ \tilde{\mathbf{V}}$}} \underbrace{\left[\begin{array}{r}
         0.737\\
        -0.675
    \end{array}\right]}_{\textstyle \text{\small $\mathbf{Q}_r$}} \hat{s}(t)
    +
    \underbrace{\left[\begin{array}{rr}
        0.737 & -0.676\\ 
        0.676 &  0.737
    \end{array}\right]}_{\textstyle \text{\small $\tilde{\mathbf{V}}$}} 
    \underbrace{\left[\begin{array}{r}
        0.675 \\
        0.737
    \end{array}\right]}_{\textstyle \text{\small $\mathbf{Q}_q$}} \hat{s}(t)^2.
\end{equation}
Interestingly, the resulting product $\tilde{\mathbf{V}}[~\mathbf{Q}_r~~\mathbf{Q}_q~]$ approximately recovers the $2 \times 2$ identity matrix. \Cref{fig:2D_R} illustrates optimized basis vectors along with the final approximation, which provides a marked improvement over any approximation in the POD basis. This gain in accuracy stems from the alignment of the secondary axis with the quadratic correction terms, enabling the model to capture significantly more nonlinear information from the original dataset.

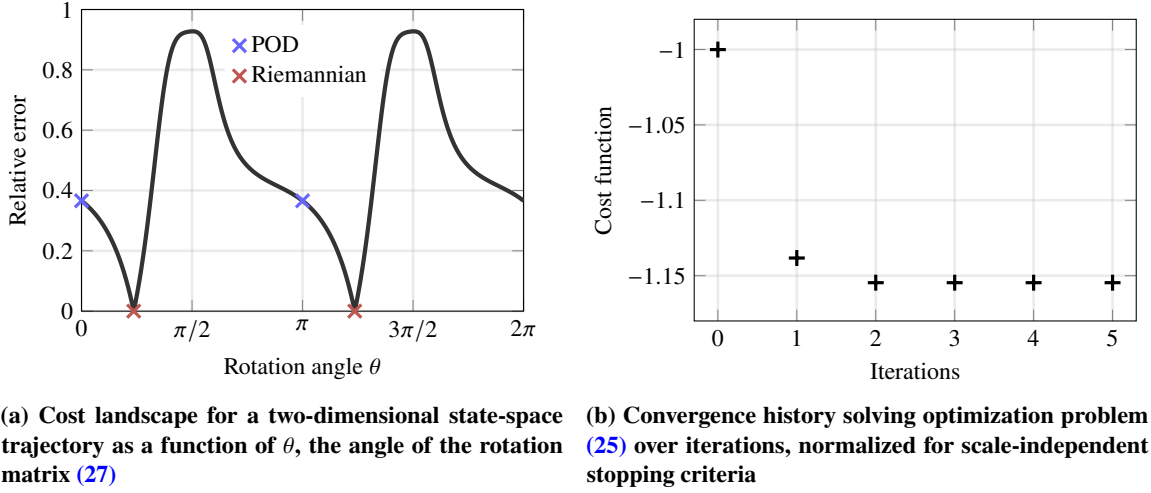
\begin{figure}[tbp]
    \centering \small
    \begin{subfigure}[t]{.45\linewidth}
    \begin{tikzpicture}
    \begin{axis}[
        width=\linewidth,
        height=.75\linewidth,
        xlabel={Rotation angle $\theta$},
        ylabel={Relative error},
        xmin=0, xmax=2*pi,
        ymin=0, ymax=1,
        mark size=1pt,
        grid=both,   
        major grid style={line width=1pt, opacity=0.3},
        xtick={0, pi/2,..., 2*pi}, 
        ytick={0, 0.2, 0.4, 0.6, 0.8, 1},
        xticklabels={$0$, $\pi/2$, $\pi$, $3\pi/2$, $2\pi$},
        legend style={
        draw=none,
        legend cell align=left, at={(0.5,.95)},anchor=north}
        ]
    \addplot[
        ultra thick,
        solid, 
        draw=black!80,
        forget plot
    ] table [
        header=false,
        x index=2,
        y index=1,
        col sep=comma
    ]{Data/2D_Example_Error_vs_Angle.csv};

    \addplot[
        only marks,
        very thick,
        mark=x, 
        draw=blue!60,
        mark size=3.5pt,
        skip coords between index={1}{360},
        skip coords between index={361}{720}
    ] table [
        header=false,
        x index=2,
        y index=1,
        col sep=comma
    ]{Data/2D_Example_Error_vs_Angle.csv};
    \addlegendentry{POD}

    \addplot[
        only marks,
        very thick,
         mark=x, 
        draw=Maroon!80,
        mark size=3.5pt,
        skip coords between index={0}{85},
        skip coords between index={86}{445},
        skip coords between index={446}{1000}
    ] table [
        header=false,
        x index=2,
        y index=1,
        col sep=comma
    ]{Data/2D_Example_Error_vs_Angle.csv};
    \addlegendentry{Riemannian}
    \end{axis}
    \end{tikzpicture}
    \caption{Cost landscape for a two-dimensional state-space trajectory as a function of $\theta$, the angle of the rotation matrix \cref{eq:2d_rotation_matrix}}
    \label{fig:Error_v_angle}        
    \end{subfigure}
    \hspace{0.5em}
    \begin{subfigure}[t]{.45\linewidth}
    \begin{tikzpicture}
    \begin{axis}[
        width=\linewidth,
        height=.75\linewidth,
        mark repeat=1, 
        xlabel={Iterations},
        ylabel={Cost function},
        xmin=-0.3, xmax=5.3,
        ymin=-1.18, ymax=-0.98,
        legend pos=north west,
        grid=both,   
        major grid style={line width=1pt, opacity=0.3},
        minor grid style={line width=.1pt, opacity=0.3}
        ]
    \addplot[
        very thick, 
        only marks,
        mark=+,
        draw=black,
        mark options={fill=black},
        mark size=3pt,
    ] table [
        header=false,
        x index=0,
        y index=1,
        col sep=comma
    ]{Data/2D_Example_Loss_vs_Iteration.csv};
    \end{axis}
    \end{tikzpicture}
    \caption{Convergence history solving optimization problem \Cref{eq:qm_least_squares5} over iterations, normalized for scale-independent stopping criteria}
    \label{fig:motivational_error}        
    \end{subfigure}
    \caption{Cost landscape and convergence history for the motivational example.}
\end{figure}

The matrix $[~\mathbf{Q}_r~~\mathbf{Q}_q~]$ can be interpreted as a rotation matrix, parameterized by an angle $\theta \in [0, 2\pi]$, applied to the column space of the candidate modes. In two dimensions, rotation matrices are parametrized as
\begin{equation}
    [~\mathbf{Q}_r~~\mathbf{Q}_q~] =
    \left[\begin{array}{rr}
        \cos(\theta) & - \sin(\theta) \\
        \sin(\theta) & \phantom{-} \cos(\theta)
    \end{array}\right],
    \label{eq:2d_rotation_matrix}
\end{equation}
where $\theta$ is the angle of rotation.
By computing the relative state error for different values of $\theta$, we can visualize how subspace orientation shapes the optimization landscape. \Cref{fig:Error_v_angle} illustrates the impact of this rotation angle on the accuracy of the quadratic approximation. Here, the optimal rotation occurs at $\theta = 0.74$ and, because the error as a function of $\theta$ is periodic with a period of $\pi$, at $\theta = 0.74 + \pi$. The optimization process rotates the basis away from a zero rotation (the original POD basis) to minimize relative error, continuing until the secondary axis aligns with the $y$-axis. At $\theta = 0.74$, the secondary axis captures the maximum nonlinear information, reducing the relative error to effectively zero. If the rotation were to continue, the next optimum would be reached at $\theta = 0.74 + \pi$; at this point, the secondary axis again aligns with the $y$-axis but in the opposite direction. Under this configuration, the nonlinear correction terms simply change sign to yield an identical approximation with the same near-zero error. \Cref{fig:motivational_error} illustrates the cost function decay using \texttt{Pymanopt}'s conjugate gradient solver under its default stopping criteria. The solver converges in six iterations upon reaching the gradient norm criterion value of $1 \times 10^{-4}$. While this motivational example is simple, the principle of optimizing subspace rotations generalizes directly to high-dimensional problems, as discussed next.

\section{Numerical Experiments on Turbulent Airfoil Wake LES}
\label{sec:numerical_experiments}

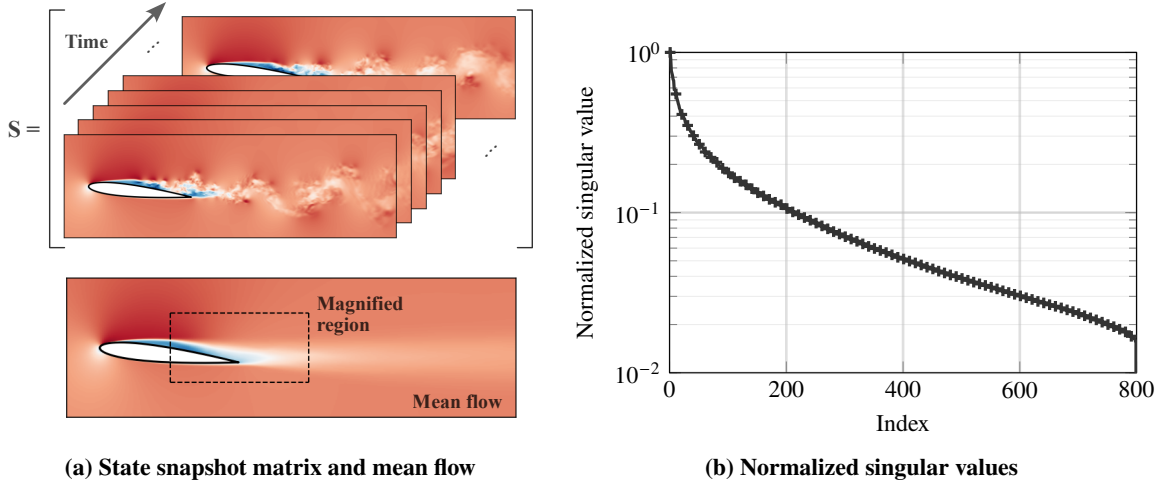
\begin{figure}[!tbp]
\centering \small
    \begin{subfigure}[t]{.42\linewidth}
        \raisebox{0.5em}{\includegraphics[width=\linewidth]{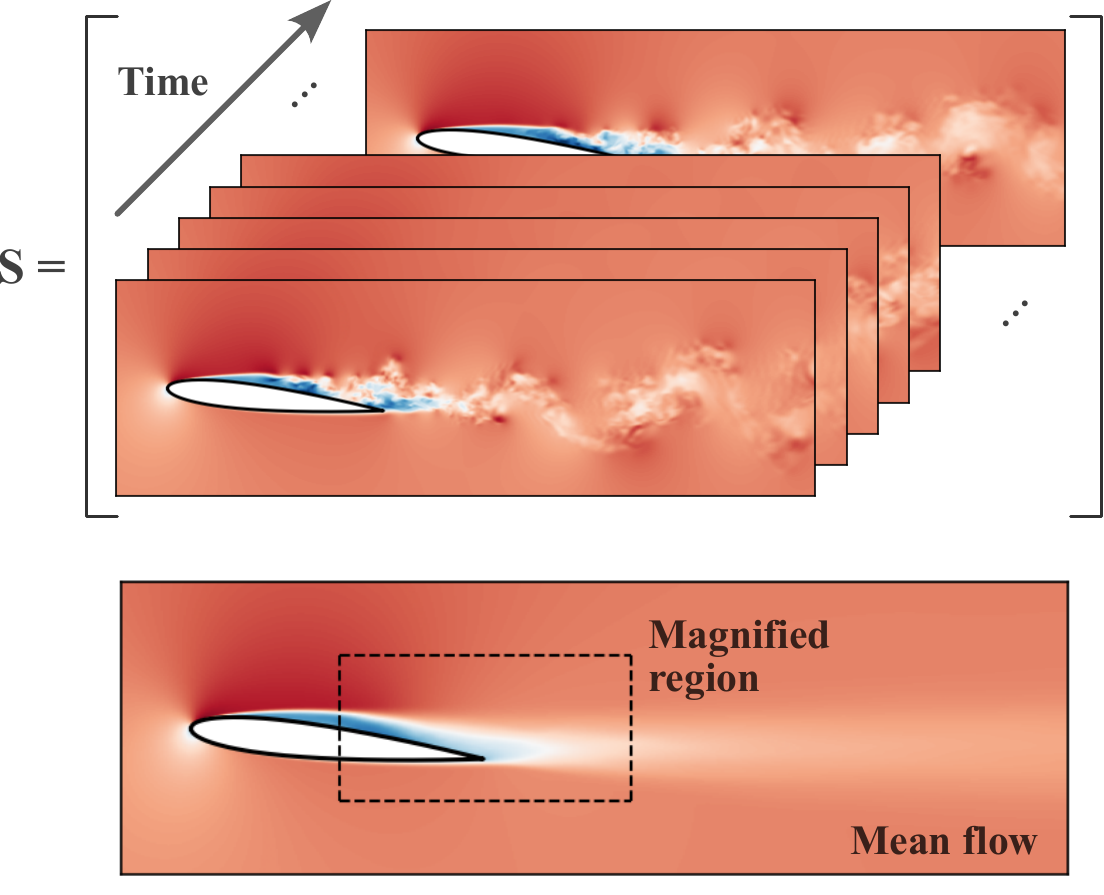}}
        \caption{State snapshot matrix and mean flow}
        \label{fig:dataset_overview}
    \end{subfigure}
    \hspace{1em}
\begin{subfigure}[t]{.47\linewidth}
    \centering \small
    \begin{tikzpicture}
    \begin{semilogyaxis}[
        width=\linewidth,
        height=.75\linewidth,
        mark repeat=10, 
        xlabel={Index},
        ylabel={Normalized singular value},
        xmin=0, xmax=800,
        ymin=1e-2, ymax=1,
        legend pos=north west,
        grid=both,   
        major grid style={line width=1pt, opacity=0.5},
        minor grid style={line width=.1pt, opacity=0.3},
        ]
    \addplot+[
        very thick,
        color=black!80,
        mark=+
        ] table [
        header=false,
        x index=0,
        y index=1,
        col sep=comma
    ]{Data/Normalized_Singular_Values_vs_Index.csv};
    \end{semilogyaxis}
    \end{tikzpicture}
    \caption{Normalized singular values}
    \label{fig:Singular_Values_LES} 
\end{subfigure}
\caption{Overview of the turbulent airfoil wake LES dataset from~\cite{deepblue}. Each snapshot provides all three velocity components extracted over the midspan slice. We plot the $x$-component of the velocity fields (left) and the corresponding singular value decay of the snapshot data (right).}
\end{figure}

To meaningfully compare the reconstruction fidelity of the different QM methods in challenging problem settings, we extract high-dimensional, time-resolved LES snapshots from the Michigan reduced-complexity modeling of fluid flows database~\cite{towne2023database}. Among the available datasets, we focus on the NACA 0012 airfoil case, which features a Mach~0.3 turbulent separated flow at a $6^\circ$ angle of attack. The Reynolds number of the flow is $23{,}000$, featuring a transitional boundary layer, separation over a recirculation bubble, and a turbulent wake. The full state dimension of the LES data is $N=666{,}630$ with $222{,}210$ grid points and the three components of the velocity field that have been non-dimensionalized by the free stream flow velocity. Rather than learning low-rank approximations for the entire temporal domain, we confine our analysis to the initial 10\% of the time window. For training and test data set construction, we use a total of $1{,}600$ snapshots along the mid-span slice. The training data set consists of 800 snapshots $\mathbf{S}_\text{train} = [~\mathbf{s}_1~~\mathbf{s}_3~~\cdots~~\mathbf{s}_{1{,}599}~] \in \mathbb{R}^{666{,}630 \times 800}$ to compute basis matrix (see \Cref{fig:dataset_overview} for a representative sample). In addition, 800 snapshots were used to construct the test data set $\mathbf{S}_\text{test} = [~\mathbf{s}_2~~\mathbf{s}_4~~\cdots~~\mathbf{s}_{1{,}600}~] \in \mathbb{R}^{666{,}630 \times 800}$. It should be noted that these training and test data sets are centered by the mean flow. The singular values of the training data are displayed in \Cref{fig:Singular_Values_LES}.

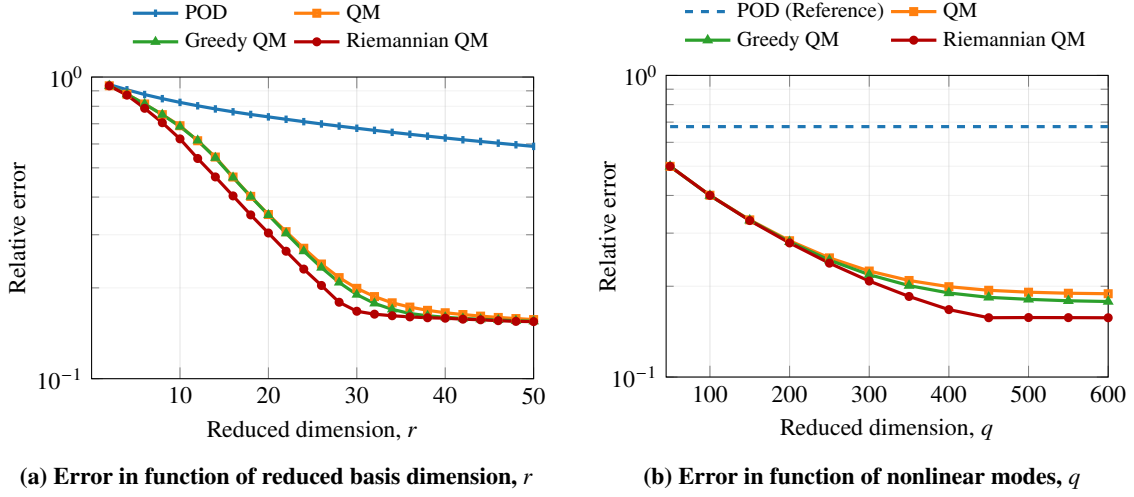
\begin{figure}[!tbp]
\centering \small
\begin{subfigure}[t]{0.45\textwidth}
    \centering
    \begin{tikzpicture}
    \begin{semilogyaxis}[
        width=\linewidth,
        height=.75\linewidth,
        xlabel={Reduced dimension, $r$},
        ylabel={Relative error},
        ylabel style={yshift=-2pt},
        xmin=0, xmax=50,
        ymin=0.1, ymax=1,
        xtick={10,20,30,40,50},
        legend style={
            draw=none,
            legend cell align=left,
            at={(0.5, 1.05)},
            anchor=south,
            font=\footnotesize,
            legend columns=2
        }, 
        grid=both,
        major grid style={line width=0.4pt, opacity=0.4},
        minor grid style={line width=0.1pt, opacity=0.2},
    ]
    \addplot[
        color=mplblue,
        mark=|,
        mark size=1.4,
        line width=1.2pt
    ] table[
        header=true,
        x index=0,
        y index=1,
        col sep=comma
    ]{Data/reduced_dimension_comparison.csv};
    \addlegendentry{POD}

    \addplot+[
        color=mplorange,
        mark=square,
        mark size=1.1,
        line width=1.2pt
    ] table[
        header=true,
        x index=0,
        y index=2,
        col sep=comma
    ]{Data/reduced_dimension_comparison.csv};
    \addlegendentry{QM}

    \addplot+[
        color=mplgreen,
        mark=triangle,
        mark size=1.1,
        line width=1.2pt
    ] table[
        header=true,
        x index=0,
        y index=3,
        col sep=comma
    ]{Data/reduced_dimension_comparison.csv};
    \addlegendentry{Greedy QM}

    \addplot+[
        color=red!70!black,
        mark=*,
        mark size=1.1,
        line width=1.2pt
    ] table[
        header=true,
        x index=0,
        y index=4,
        col sep=comma
    ]{Data/reduced_dimension_comparison.csv};
    \addlegendentry{Riemannian QM}
    \end{semilogyaxis}
    \end{tikzpicture}
    \caption{Error in function of reduced basis dimension, $r$}
    \label{fig:reduced_dimension_comp}
\end{subfigure}
\hspace{0.5em}
\begin{subfigure}[t]{0.45\textwidth}
    \centering
    \begin{tikzpicture}
    \begin{semilogyaxis}[
        width=\linewidth,
        height=.75\linewidth,
        xlabel={Reduced dimension, $q$},
        ylabel={Relative error},
        ylabel style={yshift=-8pt},
        xmin=45, xmax=600,
        ymin=0.1, ymax=1,
        legend style={
            draw=none,
            legend cell align=left,
            at={(0.5, 1.05)},
            anchor=south,
            font=\footnotesize,
            legend columns=2
        },
        grid=both,
        major grid style={line width=0.4pt, opacity=0.4},
        minor grid style={line width=0.1pt, opacity=0.2},
    ]
    \addplot[dashed,
        color=mplblue,
        mark=none,
        mark size=1.1,
        line width=1.2pt
    ] table[
        header=true,
        x index=0,
        y index=1,
        col sep=comma
    ]{Data/nonlinear_dimension_comparison.csv};
    \addlegendentry{POD (Reference)}
    \addplot+[
        color=mplorange,
        mark=square,
        mark size=1.1,
        line width=1.2pt
    ] table[
        header=true,
        x index=0,
        y index=2,
        col sep=comma
    ]{Data/nonlinear_dimension_comparison.csv};
    \addlegendentry{QM}

    \addplot+[
        color=mplgreen,
        mark=triangle,
        mark size=1.2,
        line width=1.2pt
    ] table[
        header=true,
        x index=0,
        y index=3,
        col sep=comma
    ]{Data/nonlinear_dimension_comparison.csv};
    \addlegendentry{Greedy QM}
    \addplot+[
        color=red!70!black,
        mark=*,
        mark size=1.1,
        line width=1.2pt
    ] table[
        header=true,
        x index=0,
        y index=4,
        col sep=comma
    ]{Data/nonlinear_dimension_comparison.csv};
    \addlegendentry{Riemannian QM}
    \end{semilogyaxis}
    \end{tikzpicture}
    \caption{Error in function of nonlinear modes, $q$}
\label{fig:nonlinear_dimension_comp}
\end{subfigure}

\caption{Relative test error comparison of the Riemannian QM with POD, POD-based QM, and greedy QM with respect to the number of modes: fixed $q = 400$ with varying $r$ (a), and fixed $r = 30$ with varying $q$ (b) with the fixed number of candidate modes $m = 800$.}
\label{fig:relative_err_r}
\end{figure}

To minimize the objective function on the Stiefel manifold, we utilize \texttt{Pymanopt}'s built-in conjugate gradient solver~\cite{townsend2016pymanopt}.
The optimization is governed by a termination criterion based on the norm of the Riemannian gradient; we set a convergence tolerance of $2 \times 10^{-4}$, which balances the precision of the subspace alignment with computational efficiency. To quantify the fidelity of the high-dimensional reconstructions, we use the relative state error metric
\begin{equation}
    \dfrac{\| \mathbf{S}_\text{test} - \mathbf{S}_\text{approx} \|_F} {\| \mathbf{S}_\text{test} \|_F},
\label{re_error_eq}
\end{equation}
where $\mathbf{S}_\text{test}$ represents the out-of-sample testing snapshots and $\mathbf{S}_\text{approx}$ denotes the approximations generated by means of low-rank approximation through POD~(\Cref{subsec:pod}), POD-based QMs (\Cref{subsec:pod_based_qm}), Greedily constructed QMs (\Cref{subsec:greedy_qm}), and Riemannian QMs from \texttt{FastQM} (\Cref{sec:riemann}).

\Cref{fig:reduced_dimension_comp} shows the relative state error as a function of the reduced basis dimension $r$. To mitigate overfitting, the regularization parameter is set to $\gamma = 10^2$ for all quadratic approximations. The number of nonlinear modes is fixed at $q = 400$, and the number of candidate modes at $m = 800$. Because standard POD does not include correction terms, it serves as a baseline for evaluating the performance of the quadratic manifold approximations. Note that the slow decay of the singular values makes this a particularly challenging test case. We start with the observation that all quadratic manifold approaches outperform the standard POD baseline. It can also be seen that the greedy QM delivers nearly identical performance to the POD-based QM at low dimensions; this is because its first eight selected indices match the leading singular vectors exactly, neutralizing any advantage provided by greedy selection. The proposed Riemannian QM approach, on the other hand, avoids this limitation through optimization on the Stiefel manifold, generating basis vectors that are not strictly bound to the dataset's left singular vectors. As a result, it consistently achieves the highest accuracy across the full range of $r$-values. \Cref{fig:nonlinear_dimension_comp} details the impact of varying the number of nonlinear modes, $q$, at a fixed dimensionality of $r=30$. While incorporating more modes initially improves the quadratic correction, accuracy plateaus around $q=400$. Mirroring the fixed-$q$ trend, the Riemannian approach maintains the lowest error throughout the parameter sweep.

\begin{table}[!t]
\centering \small
\caption{Impact of number of candidate modes, $m$, and the number of quadratic correction terms, $q$, on relative test error ($r=30$).}
\label{tab:candidate_modes_comp}
\setlength{\tabcolsep}{10pt} 
\renewcommand{\arraystretch}{1.1} 
\begin{tabular}{cccccc}
\toprule
Candidate modes, $m$ & $q=250$ & $q=350$ & $q=450$ & $q=550$ & $q=600$ \\
\midrule
300 &  0.2384  &  -       &  -      &  -    &    -    \\
400 &  0.2382  &  0.1843  &  -      &  -    &    -    \\
500 &  0.2381  &  0.1840  &  0.1580  &  -    &    -    \\
600 &  0.2381  &  0.1837  &  0.1567  &  0.1559  &  - \\
700 &  0.2382  &  0.1839  &  0.1563  &  0.1550  &  0.1550  \\
800 &  0.2382  &  0.1840  &  0.1560  &  0.1552  &  0.1547  \\
\bottomrule
\end{tabular}
\end{table}

\Cref{tab:candidate_modes_comp} investigates the impact of the number of candidate/feature modes, $m$, on approximation accuracy. For this assessment, the reduced dimension is fixed at $r = 30$ across various choices for the number of nonlinear modes, $q$. Some data points are unavailable at certain values of $m$; this omission occurs because the number of candidate modes $m$ must exceed the sum of the linear and nonlinear modes, $r + q$, to satisfy the linear independence condition stated in \cref{eq:linear_independence}. Where this condition is met, the relative error remains nearly constant. This demonstrates that restricting the search space to a predefined feature span of the leading $m$ left singular vectors of the data matrix is a logical choice, provided that $m$ is sufficiently large but not overly conservative.

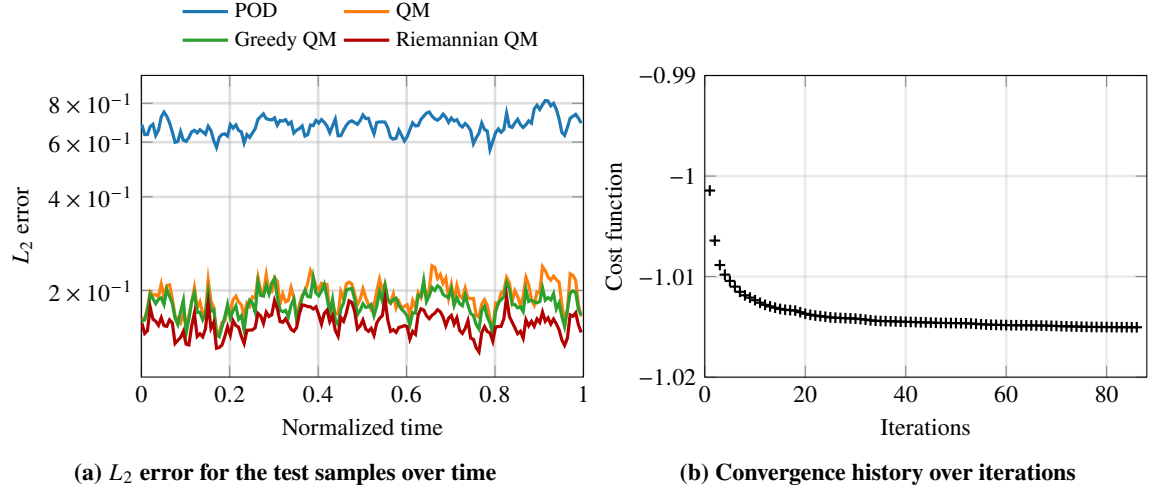
\begin{figure}[!tbp]
    \centering \small
    \begin{subfigure}[t]{.45\linewidth}
    \begin{tikzpicture}
    \begin{semilogyaxis}[
        width = \linewidth,
        height = 0.75\linewidth,
        xlabel={Normalized time},
        ylabel={$L_2$ error},
        xmin=0, xmax=1,
        mark size=1pt,
        legend style={
            draw=none,
            legend cell align=left,
            at={(0.5, 1.05)},   
            anchor=south,       
            font=\footnotesize,
            legend columns=2   
        },        
        grid=both,   
        major grid style={line width=1pt, opacity=0.5},
        minor grid style={line width=.1pt, opacity=0.3},
        ytick={0.2, 0.4, 0.6, 0.8},
        yticklabels={$2 \times 10^{-1}$, $4 \times 10^{-1}$, $6 \times 10^{-1}$, $8 \times 10^{-1}$},
        yticklabel style = {font=\small}
        ]
    \addplot[
        mplblue,
        no marks,
        line width=1.2pt,
        each nth point={5}
    ] table [
        header=true,
        x expr={\thisrowno{0}/1600},
        y index=1,
        col sep=comma
    ]{Data/time_data_500.csv};
    \addlegendentry{POD}

    \addplot+[
        mplorange,
        no marks, 
        line width=1.2pt,
        each nth point={5}
    ] table [
        header=true,
        x expr={\thisrowno{0}/1600},
        y index=2,
        col sep=comma
    ]{Data/time_data_500.csv};
    \addlegendentry{QM}

    \addplot+[
        mplgreen,	
        no marks, 
        line width=1.2pt,
        each nth point={5}
    ] table [
        header=true,
        x expr={\thisrowno{0}/1600},
        y index=3,
        col sep=comma
    ]{Data/time_data_500.csv};
    \addlegendentry{Greedy QM}
    
    \addplot+[
        color=red!70!black,	
        no marks, 
        line width=1.2pt,
        each nth point={5}
    ] table [
        header=true,
        x expr={\thisrowno{0}/1600},
        y index=4,
        col sep=comma
    ]{Data/time_data_500.csv};
    \addlegendentry{Riemannian QM}
    \end{semilogyaxis}
    \end{tikzpicture}
    \caption{$L_2$ error for the test samples over time}
    \label{fig:L2_err_time}
    \end{subfigure}
    \hspace{0.8em}
    \begin{subfigure}[t]{.45\linewidth}
    \begin{tikzpicture}
    \begin{axis}[
        width=\linewidth,
        height=.75\linewidth,
        mark repeat=1, 
        xlabel={Iterations},
        ylabel={Cost function},
        xmin=0, xmax=88,
        ymin=-1.02, ymax=-0.99,
        grid=both,   
        major grid style={line width=1pt, opacity=0.3},
        yticklabel style = {font=\small}
        ] 
        \addplot[
        only marks,
        color=black,
        mark=+,
        mark options={fill=black},
        line width=0.8pt
        ] table[
        col sep=comma,
        header=true,
        x index=0,
        y index=1,
        ]{Data/cost_func3.csv};
    \end{axis}
    \end{tikzpicture}
    \caption{Convergence history over iterations}
    \label{fig:iteration_error_les}        
    \end{subfigure}
    \caption{Temporal evolution of the $L_2$ error and convergence behavior for the proposed optimization approach at $r = 30$. The approximation uses $q = 500$ optimized modes for the quadratic term $\Vq$, optimized from from $m = 800$ candidate left singular vectors.}
    \label{fig:L2error_convergence}
\end{figure}

\Cref{fig:L2_err_time} presents the $L_2$ error over time for the state reconstruction in the original state space, evaluated at $r=30$, $q=500$, and $m=800$. Consistent with our previous findings, the proposed Riemannian approach yields the lowest reconstruction error over time among all dimensionality reduction methods considered. As illustrated in \Cref{fig:iteration_error_les}, this superior accuracy is achieved in 86 iterations, requiring about 20~seconds to satisfy the gradient norm convergence criterion of $2 \times 10^{-4}$. Because the original cost function magnitude is excessively large, on the order of $10^{6}$, the objective function was normalized by its initial value. Furthermore, the objective function plateaus and ceases to show meaningful variation after a certain number of iterations, indicating that the optimization successfully reaches a well-converged solution.

\begin{figure}[!tbp]
    \centering
    \includegraphics[width=0.4\linewidth]{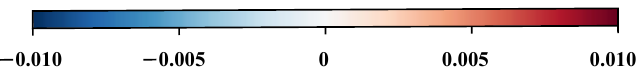} \\[0.25em]
    \begin{subfigure}[t]{0.49\textwidth}
            \includegraphics[width=\linewidth]{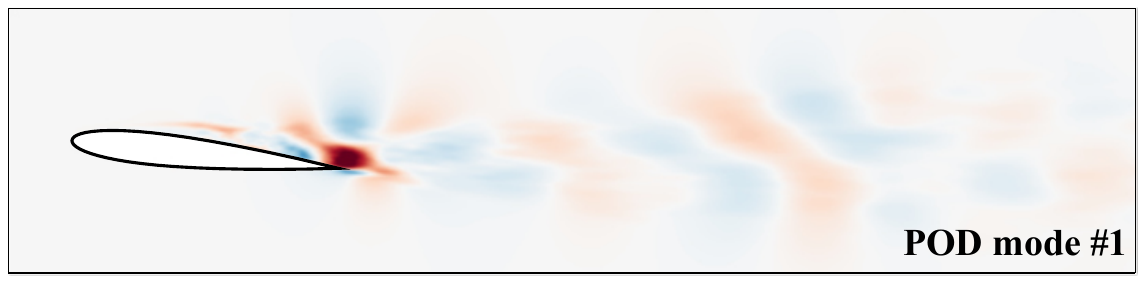} \\[0.25em]
            \includegraphics[width=\linewidth]{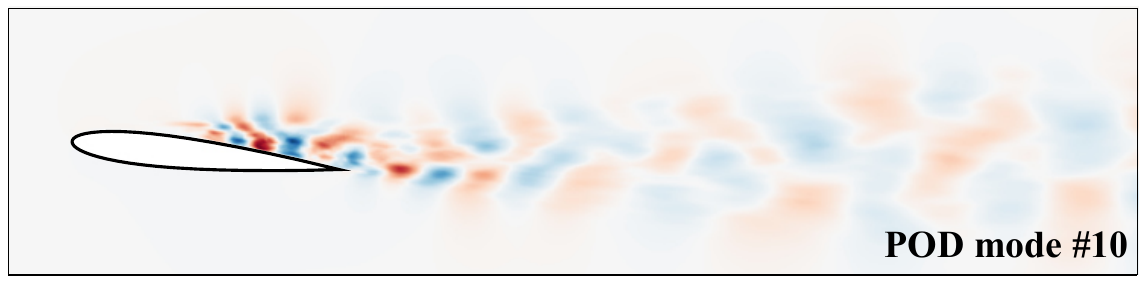}  \\[0.25em]
            \includegraphics[width=\linewidth]{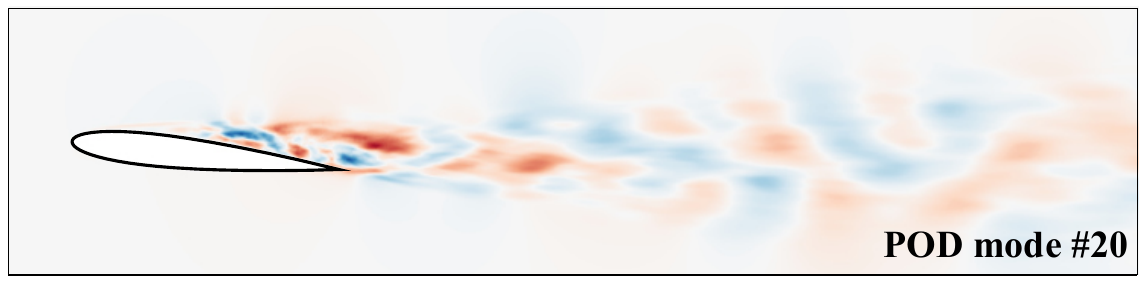}
        \caption{POD modes}
        \label{POD_Modes}
    \end{subfigure}
    \begin{subfigure}[t]{0.49\textwidth}
    \centering
            \includegraphics[width=\linewidth]{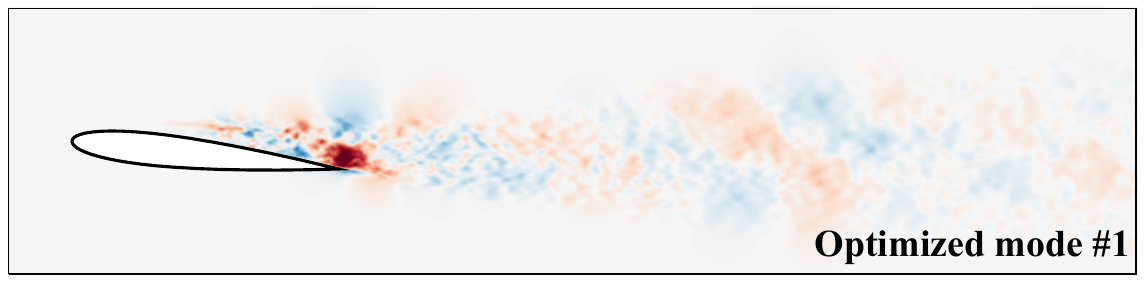}   \\[0.25em]
            \includegraphics[width=\linewidth]{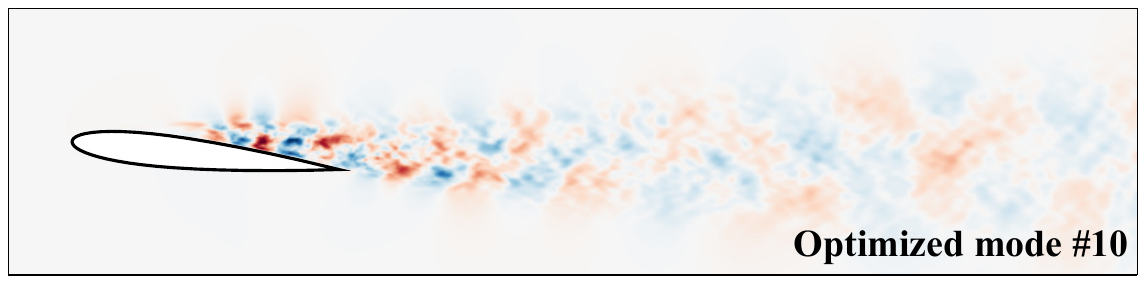}  \\[0.25em]
            \includegraphics[width=\linewidth]{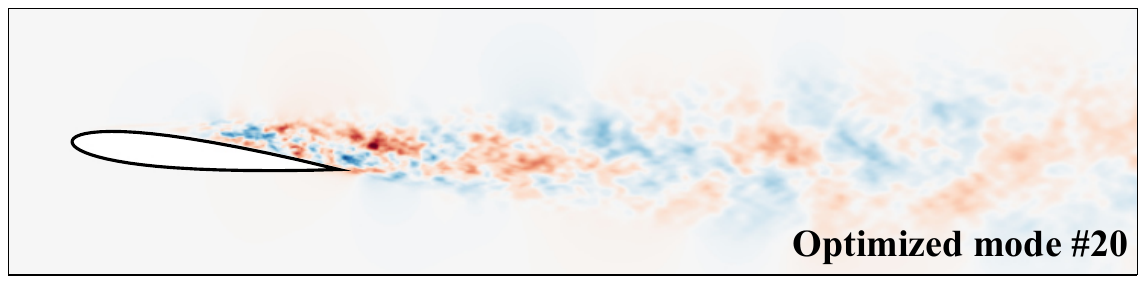}
        \caption{\emph{Optimized} modes}
        \label{Optimized_POD_Modes}
    \end{subfigure} 
    \caption{Comparison of the spatial modes at a dimensionality of $r = 30$. The number of nonlinear modes (columns of $\Vq$) is chosen to be $q = 500$ with each optimized mode represented as a linear combination of $m=800$ candidate left singular vectors.}
    \label{modes_comp} 
\end{figure}

\Cref{modes_comp} compares the leading POD modes (the columns of $\mathbf{\tilde{V}}$) against the optimized modes, the columns of $\mathbf{\tilde{V}}[~\mathbf{Q}_r~~\mathbf{Q}_q~]$. We focus on the 1st, 10th, and 20th modes. These basis vectors are constructed using $r=30$ and $q=500$ and are chosen from a candidate span of $m = 800$. Notably, while the standard POD modes in \Cref{POD_Modes} exhibit broader, fuller structures, the optimized modes in \Cref{Optimized_POD_Modes} resolve significantly finer scales in a way that facilitates quadratic approximation. This structural difference allows the Riemannian quadratic manifold to capture critical nonlinear features which are simply not present in the leading POD modes.

\begin{figure}[!tbp]
\centering
\includegraphics[width=.875\linewidth]{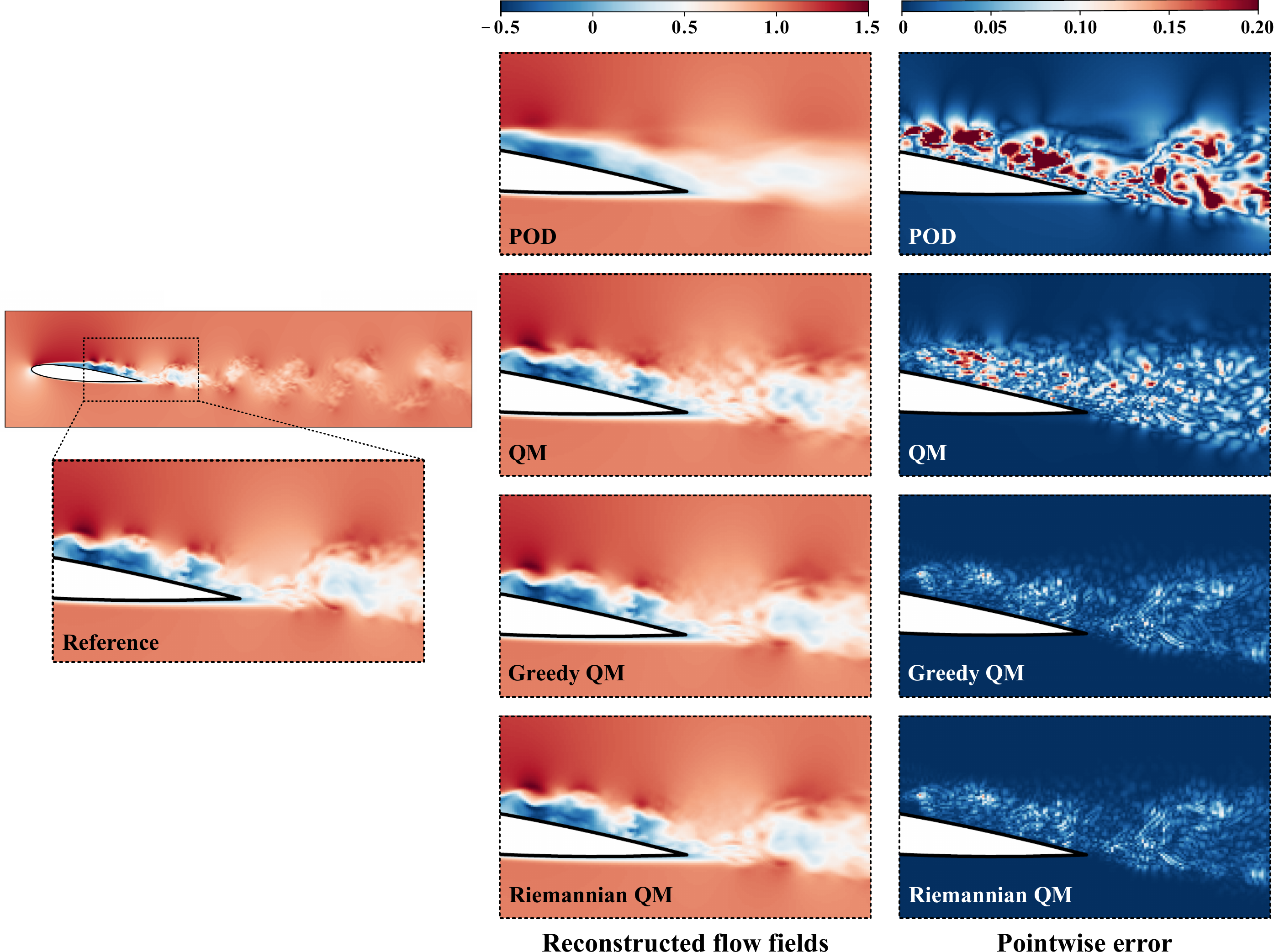}
\caption{Snapshots of the $x$-velocity component and pointwise reconstruction errors at normalized time $t=0.5$. Reconstructions (middle) and error distributions (right) for $r=30, q=500, m=800$ are compared against the reference solution (left).}
\label{fig:comparison}
\end{figure}

To visualize error across the computational domain, the reconstructed velocity field and pointwise error distribution for a test snapshot at normalized time $t=0.5$ are shown in \Cref{fig:comparison}. The linear POD baseline struggles to capture the wake structures near the airfoil's trailing edge, which is reflected in its high reconstruction error. The POD-based QM offers an improved depiction of the wake structure, yet still exhibits noticeable pointwise errors. Both the greedy QM and the proposed Riemannian QM accurately capture the unsteady wake dynamics, reducing the relative state-vector errors to 2.22\% and 2.06\%, respectively. While the proposed method shows only a slight advantage for this individual state vector, its superior overall performance and robustness are clearly demonstrated by the global results across the entire temporal domain.

\section{Concluding Remarks}
\label{sec:concluding_remarks}

This work introduces \texttt{FastQM}, a novel dimensionality reduction framework treating quadratic manifold (QM) construction as a continuous optimization problem on the Stiefel manifold. While model reduction via QMs offers a promising route toward capturing complex dynamics, existing formulations struggle to maximize the effectiveness of quadratic correction terms. This suboptimality stems from their reliance on standard POD subspaces, which can misalign with the underlying nonlinear geometric structures. \texttt{FastQM} overcomes this bottleneck by utilizing orthogonal rotation matrices computed via Riemannian optimization, enabling explicit coordinate transformations that improve quadratic approximation accuracy. The performance and robustness of the proposed framework were evaluated using a high-dimensional LES dataset of a turbulent airfoil wake. Numerical experiments show that \texttt{FastQM} consistently achieves superior reconstruction fidelity compared to its POD-based counterparts, demonstrating that the effectiveness of quadratic corrections can be optimized at a manageable computational footprint. Future work will extend this framework to a wider variety of engineering problems to further assess its generalizability and performance. 


\bibliography{bib}

\end{document}